\documentclass[12pt]{article} \usepackage{amsmath,amsthm,amssymb,color}
\usepackage{graphicx}
\usepackage{enumerate}
\usepackage{natbib}
\usepackage{url} 

\newcommand{\blind}{0}

\newcommand{\Expec}{\mathbb E}
\newcommand{\Prob}{\mathbb P}
\newcommand{\Var}{{\mathbb V}{\mathrm ar}}
\newcommand{\Cov}{{\mathbb C}{\mathrm ov}}

\begin{document}

\def\spacingset#1{\renewcommand{\baselinestretch}%
{#1}\small\normalsize} \spacingset{1}


\if0\blind
{
  \title{\bf Detection and Estimation of Local Signals}
  \author{Xiao Fang\\
    Department of Statistics, The Chinese University of Hong Kong\\
    and \\
    David Siegmund \\
    Department of Statistics, Stanford University}
  \maketitle
} \fi

\if1\blind
{
  \bigskip
  \bigskip
  \bigskip
  \begin{center}
    {\LARGE\bf Detection and Estimation of Local Signals}
\end{center}
  \medskip
} \fi

\bigskip
\begin{abstract}
We study the maximum score statistic to detect and estimate local signals in the form
of change-points in the level, slope, or other property of a sequence
of observations, and to segment the sequence when there appear to be multiple
changes. We find that when observations are serially dependent, the
change-points can lead to upwardly biased estimates of autocorrelations,
resulting in a sometimes serious loss of power.
Examples involving temperature variations, the level of 
atmospheric greenhouse gases,
suicide rates, incidence of COVID-19, and excess deaths during the pandemic
illustrate the general theory.
\end{abstract}

\noindent%
\vfill

\newpage
\spacingset{1.5} 

\section{Introduction.} 
We consider a problem of detection, estimation, and segmentation 
of local nonlinear signals imbedded in a sequence of observations.  
As a  model, we assume for $u = 1, \ldots, T$
\begin{equation} \label{loglik1}
Y_u = \rho Y_{u-1} +\mu(X_u) +  \sum_k\xi_k f[(X_u-t_k)/T] + \epsilon_u,
\end{equation}
where $t_1 < t_2 < \ldots < t_M$  define the locations and $f$
the ``shape,'' 
of the local signals. For asymptotic analysis given below,
we assume that the $t_k$ are scaled, so $t_k/T = t_{0,k}$.  
Initially we assume the $X_u$ are fixed,
but for some applications they are random.   The $\epsilon_u$ are
independent mean 0, normally distributed errors.  
For the moment we assume their variances are known and equal one, and
discuss later how they should be estimated.
The nuisance parameters $\mu(X_u)$
depend on the variable $X_u$ and may be
constant or a parameterized regression function.
The purpose of introducing $\rho$ into our model is to decrease the
risk of false positive errors if the observations are dependent, but
it is not intended to be a good description of the true dependence.
For our general theory we assume $\rho$ is an unknown
constant, but for some applications we regard $\rho$ as a 
pre-whitening device and assign it an arbitrary numerical value, which
may be altered in
subsequent analysis if we appear to have made an inappropriate choice. 

An important special case that we return to in examples below
is $X_u = u$ and $\mu_u = \alpha +\beta[(u-(T+1)/2)/T]$, 
so with this notation the model becomes 
\begin{equation} \label{loglik}
Y_u = \rho Y_{u-1} +\mu_u +  \sum_k\xi_k f[(u-t_k)/T] + \epsilon_u.
\end{equation}
It is sometimes convenient to simplify several basic  
calculations by considering an alternative continuous time  
model, for which (\ref{loglik}) can be conveniently written
\begin{equation} \label{cont}
d Y_u = -\gamma Y_u du + \mu_u du +  \sum_k\xi_k f[(u-t_k)/T]du + dW(u) 
\end{equation}
for $0 \leq u \leq T$,
where $dW$ defines white noise residuals.  

A feature of all these models is statistical irregularity: $\xi_k$ and
$t_k$ are confounded in the sense that if $\xi_k = 0$, then $t_k$ is
not defined.
Specific choices of the function $f$ appear to be appropriate for a
variety of applications, some of which are discussed below.  
The special case $f(x) = I\{x > 0\}$ or $I\{0 < x < \tau\}$ is a 
frequently discussed
``change-point'' model, which has been
applied to a variety of problems, usually with $\rho$ assumed equal to 0.
See, for example, \cite{OlVe04}, \cite{Lund16} and references given there.  
\cite{FaLSi20} and \cite{Fr14} address a version of the problem that
involves the possibility of  multiple change-points, and a major goal is segmentation of the 
observations to identify their number and locations, while controlling the probability
of false positive detections.  In this paper we 
extend the methods of \cite{FaLSi20} to deal with more general signals 
and with
observations that may have  autoregressive dependence.
See also \cite{BCFr16}.

    
A model with a long history, primarily for independent 
observations and at most one break-point, is a broken line regression model,
where $f(x) = x^{+}$, so the model is given by
\begin{equation} \label{broken}
Y_u =\rho Y_{u-1} + \alpha + \beta [(u-(T+1)/2)/T] + \sum \xi_k [(u-t_k)^+/T] + \epsilon_u.
\end{equation}
An application to monitor kidney function after 
transplant has been
discussed in a an elegant series of papers 
by A. F. M.  Smith and others (e.g. \cite{SC80}). 
\cite{Dav87} and \cite{KnSi89} give theoretical analyses for 
independent observations and at most one change-point.
\cite{ToLes03} provides an analysis and ecological applications.
Several of the examples discussed below involve climatological 
time series, or day by day
newly confirmed cases of COVID-19.
Although alternative models are possible, 
a broken line model provides a simple conceptual framework that 
allows us to determine
whether putative changes in direction of a regression 
function are real, and to  
suggest or confirm hypotheses of scientific interest.  In particular,
we can evaluate the evocative ``hockey stick'' suggested by the 
appearance of  
annual average temperatures of many countries in the 20th century.

In a variety of genomic applications
$t$ denotes a genomic location, and different applications suggest
functions $f$  having different characteristic shapes.  For example, 
for detection of copy number variations (CNV) the indicator of
an interval or a half line may be appropriate (e.g., 
\cite{OlVe04}, \cite{FaLSi20}).
For  detection of 
differentially methylated genomic regions
\cite{Irizarry12} suggests a ``bump'' function.
See simple examples in Appendix D of the Supplementary Material  and 
elaborations of this model 
involving covariates in \cite{JaffeIrizarry12}.
For ChIP-Seq analysis steeply decaying 
symmetric functions like $(1-|u|)^+$,  $\exp(-|u|)$,
or a normal probability density function are different possibilities
(cf. \cite{Liu13}, \cite{SJGM2013}). 
In some applications it is often appropriate to add a scale parameter
$\tau$, so at $u$ the signal located at $t$ is of the form $\xi f[(u-t)/(\tau T)]$
(cf. \cite{SiWo95}).  The model of ``paired change-points''
of \cite{OlVe04} can be described similarly, with $f(x) = I\{0 < x \leq 1\}$.

In the special case that $\mu$ is a constant, $X_u = Y_{u-1}$, and $f(x)$ is 
the indicator that
$x \leq  0,$ (\ref{loglik1}) is the simplest example of a 
threshold autoregression,
as introduced by Tong and studied by Chan and Tong and others 
(e.g., \cite{CT90} ) in numerous projects.  
In a still different context,
$X_u$ may be a regressor, say a biomarker,
and the local signal can be used to study whether a subset of individuals, 
defined by the
value of that biomarker, differ from others in some respect, e.g.,
response to a treatment or disease susceptibility.   

In some applications $X_u$ may be multidimensional
or the noise distribution may not be normal, with the Poisson distribution 
representing a particularly interesting alternative.  

In Section 2 we begin 
with some basic calculations and a discussion of testing
for at most one change. We find two important
features of our method.  (i) The maximum likelihood estimator of $\rho$
under the hypothesis of no local signals can be badly biased when there are 
signals. (ii) The variance of the score statistic under the hypothesis of
no signal, hence asymptotic evaluation of false positive error 
probabilities when $\rho$
is estimated by maximum likelihood does not depend on the value of $\rho$.
We also develop
confidence regions for $t$, which add appreciably to our 
understanding, and joint regions for $t$ and $\xi$. 

A novel theoretical contribution are methods
of segmentation, 
discussed in Section 3 and illustrated in 3.2. 
Two of the three methods suggested there are 
conceptually similar to methods proposed in \cite{FaLSi20} for 
jumpchanges with independent observations, but their 
implementation is based on new probability approximations that 
involve random fields having smooth behavior in one coordinate 
and not smooth behavior
in others.
The third is a new
procedure, which like circular binary segmentation  
(CBS,\cite{OlVe04}) for the detection of jump changes,
is designed to detect pairs of changes, but unlike CBS does not 
assume that the paired changes move in opposite directions.

If we assume as a consequence of segmentation that the $t_k$ 
in (2) are known without error,
the model becomes a linear model.  An important ingredient
of our overall strategy is to use that linear model 
to evaluate the global success of our segmentation and to re-estimate
the value of $\rho$, which is a lagged variable in
the regression function.  This provides the possibility of a 
re-segmentation if initially we used
an inappropriate value.  
See Section 3.1 for examples.

Section 4 briefly
considers the use of our methods for sequential detection of a 
change of slope in broken line
regression, with examples from the COVID-19 pandemic.  

Additional theoretical calculations and examples are 
contained in Supplementary Material.

Much although not all the related literature assumes independence, which 
seems adequate for many applications.  We use a very simple model of dependence,
defined in terms of the observed process, which makes that process first order
Markov.  Our goal here is not to build an 
approximately correct dependence structure, but
to protect our methods for detecting changes in the mean value against an 
excess of false positives when the
process exhibits dependence.  At the cost of some 
computational complications we can also
consider second or higher order Markov dependence. 

There are alternative models of dependence, usually defined in terms of
the unobserved residuals.
See
\cite{Lund16} for an interesting discussion 
using weak convergence arguments that involve 
rescaling the temporal and spatial coordinates. 
Explicit short term dependencies vanish in the limit, because points
in the time scale that are separated by $o(m)$ are 
merged into a single point.  This makes interpretation 
difficult if there appear to be nearby change-points in 
the original time scale.  

Still other  methods focus on estimation of the
regression function (e.g., by a penalized likelihood), but do not
provide statistical measures of the validity and variability 
in location of detected signals.

Our simple Markov model seems to provide a flexible 
technique to protect against an excess of false positive
errors without a substantial loss of power.  


\section{Basic Calculations and the Case of at Most One Change.}
To introduce our notation and provide some basic results, we
begin with the  model given by (\ref{loglik}), with 
$M = 1$ or 0, and we consider a test
of the hypothesis $\xi = 0$.  The log likelihood is given by $\ell(\xi, t, \theta)$
\[= -\frac{1}{2} \sum_u \{ Y_u - \rho Y_{u-1} - \alpha - \beta [(u-(T+1)/2)/T] - \xi f[(u-t)/T]\}^2,
\]
where $\theta = (\alpha, \beta,\rho)'$, and we consider a test 
based on the maximum with respect to
$t$ of the standardized score statistic
\begin{equation} \label{score}
\ell_\xi (0,t,\hat{\theta})/\sigma(t),
\end{equation}
where $\ell_\xi = d \ell/d \xi$, 
$\hat{\theta}$ is the maximum likelihood estimator of
the nuisance parameters
under the hypothesis $\xi = 0$, and $\sigma^2(t)$ is 
the asymptotic variance of the numerator.

Our large sample theory begins with the standard expansion of the
numerator in (\ref{score}), given by
\begin{equation} \label{score1}
\ell_\xi (0,t,\hat{\theta}) \approx  \ell_\xi (0,t,\theta) - I_{\xi,\theta} 
I_{\theta,\theta}^{-1} \ell_\theta,
\end{equation}
where $I_{\cdot,\cdot}$  denotes elements of the Fisher information 
matrix and by the law of large numbers all 
quantities on the right hand side are evaluated at $\xi = 0$ 
and true
values of the other parameters.  This expansion is valid in 
large samples up to terms that are $o(T)$  in probability.  To 
emphasize the structure of this approximation
we use the notation
$\Expec_0$ to denote expectation under the hypothesis $\xi = 0$ and write
(\ref{score1}) in the form
$V_t =  \ell_\xi(t) - \Psi(t)'A \ell_\theta$,
where $\ell_\xi(t) = \ell_\xi (0,t,\theta)$,
$\Psi(t)' = \Expec_0[\ell_\xi(t) \ell_\theta']$ is not random although it
depends on $t$; $\ell_\theta$ does not depend on $t$, and 
under the hypothesis $\xi = 0$  it has mean value 0
and covariance matrix $A^{-1} = I_{\theta,\theta}$. 

\medskip\noindent
{\bf Remark.}  It may be shown that the decomposition given in (6)
does not depend on the normality of the $\epsilon_j$, although we
can no longer use the terminology of likelihood, efficient score, etc.,
and must rely on the central limit theorem to justify the asymptotic 
normality of probability calculations given below.  

Calculations yield a number of simple propositions.

\medskip\noindent
{\bf Proposition 1.}
\begin{equation} \label{covariance}
\Sigma(s,t)  = \Expec_0(V_s V_t) = \Expec_0 [\ell_\xi(s) \ell_\xi(t)] 
- \Psi(s)'A\Psi(t).
\end{equation}
In particular the  variance of $V_t$ is $\sigma^2(t) = \Sigma(t,t).$

\smallskip\noindent
{\bf Remarks.} (i) Additional calculations show that $\Sigma(s,t)$ does
not depend on unknown parameters.  While this is expected and easily
demonstrated for $\alpha$ and $\beta$,
we found it tedious to demonstrate (cf. the online
Supplementary Material) in the case of $\rho$. 
(ii) In the case the local signal contains a scale parameter,
so it takes the form $\xi f[(u-t)/(\tau T)]$,
the covariance is  similar in its general formulation, but slightly more complicated 
to compute in examples.

Writing $\Expec_{t,\xi}$ to denote expectation under the alternative, where $t$ 
and $\xi$ can be  vectors,
we have

\medskip\noindent
{\bf Proposition 2.}
$\Expec_{t,\xi}(V_s) = \sum_k \Sigma(s, t_k) \xi_k.$

Using Proposition 2, we can calculate numerically the expected sample path of
$V_s$ or $Z_s = V_s/\sigma_s$ to see how we might go about detecting local signals.
In the case that there is only one signal and, say, $\xi > 0$, the expected path is typically
maximized at $t$.
Hence, if we assume there is at most one change at some 
unknown value $t_0$, we put
$Z_t = [\ell_\xi(t) - \Psi(t)'A \ell_\theta]/\sigma(t)$,
and use the maximizing value of $|Z_t|$ as an estimator of
$t_0$, provided $\max_t |Z_t| > b$ for a suitable threshold $b.$
(In Section 3 we discuss the segmentation of the values of $t$ in the case of multiple signals
by considering local extrema of $|Z_t|$, judged against appropriate local backgrounds.)

For the case of at most one change, the false positive error
probability  is
$\Prob_0 \{ \max_t |Z_t| \geq b \}$.
To evaluate this probability, we distinguish essentially two cases, 
one where $f$ is continuous and the second where it is discontinuous.
The case where $f$ is the indicator of a half line or an
interval (of unspecified length) is discussed in
\cite{FaLSi20}.  Here we assume that $f$ is continuous
and piecewise differentiable. 
An approximation based on
Rice's formula along the lines suggested by 
Hotelling (1939), Davies (1976) and others is given by

\medskip\noindent
{\bf Proposition 3.}
\begin{equation} \label{Rice1}
\Prob_0\{ \max_{T_0 \leq t \leq T_1} |Z_t| > b\} \leq 2\{ (\varphi(b)/(2\pi)^{1/2}) \int_{T_0}^{T_1}
[\Expec_0(\dot{Z}_t)^2]^{1/2} dt
+ 1 - \Phi(b)\},
\end{equation}
where $\dot{Z}_t$ denotes the derivative of $Z$ at $t$.
For the model where the local signal at $u$ has the form $\xi f[(u-t)/ \tau T]$,
and we maximize over both $t$ and a range of values of $\tau$, a first order approximation 
for large $b$ is
\begin{equation} \label{Rice2}
\Prob_0\{ \max_{t,\tau} |Z_{t,\tau}| > b\} \sim 2\{ b \varphi(b)/(2\pi)\} \int_{\tau_0}^{\tau_1}  \int_{T_0}^{T_1}
\rm{det}[\Expec_0(\dot{Z}\dot{Z}')]^{1/2} dt d\tau,
\end{equation}
where $\dot{Z} = \nabla Z_{t,\tau}$.
This approximation can be improved by adding lower order terms 
involving edge effects and corrections for
curvature.  The  most important is the boundary correction 
at the minimum value of $\tau$ equal to
$ (8\pi)^{-1/2} \varphi(b) \int_{T_0}^{T_1} [\Expec_0(\dot{Z})^2]^{1/2} dt$, where
now $\dot{Z} = \frac{\partial Z}{\partial t} ( t, \tau_0)$.  It is often convenient
when $T$ is large to ignore edge effects.
See \cite{SiWo95} or \cite{AT07} for more detailed approximations. 

\smallskip
To evaluate (\ref{Rice1}) and (\ref{Rice2}), the following result is useful.

\smallskip\noindent
{\bf Proposition 4.}
\begin{equation}
\Expec_0[(\dot{Z}\dot{Z}')] = [\Expec_0(\dot{\ell} \dot{\ell}') - \dot{\Psi}'A \dot{\Psi}]/\sigma^2(t) 
- [\dot{\sigma}\dot{\sigma}']/\sigma^2. 
\end{equation} 

Finally observe that by  Propositions 1 and 2, we have the ``matched filter'' conditions,
and consequently

\smallskip\noindent
{\bf Proposition 5.}
Given $t$ (or $t,\tau$),  the variable $V_t$ (or $V_{t,\tau}$) is sufficient for $\xi$.

\smallskip
For the specific case of (\ref{broken}), it simplifies calculations 
somewhat to consider the continuous time 
version given in (\ref{cont}), which amounts to replacing certain sums by
integrals, leading to

\medskip\noindent
{\bf Proposition 6.}
$\sigma^2(t) \sim T\{(1-t_0)^3/3 - (1-t_0)^4(1+t_0 +t_0^2)/3\}$,
and $\Expec_0(\dot{V}_t^2) \sim T\{(1-t_0) -(1-t_0)^2 - 3[t_0(1-t_0)]^2\}.$

\subsection{Confidence Regions}
In this section we discuss confidence regions for the parameters
$t$ or  joint regions for $t$  and $\xi$.
Calculation shows that the expected value of
$V_s$ equals $\xi$ times the covariance of $V_s$ and $V_t$ (cf. Proposition 2).
It follows from straightforward multivariate calculation
that $Z_t = V_t/\sigma(t)$ is sufficient for $\xi$,
so the conditional distribution of $Z_s$ given $Z_t$ is the same
as it would be under the null hypothesis, $\xi = 0$.

The case of jump changes was studied by
\cite{FaLSi20}, so here we consider the case of continuous processes.
For a similar approach based on a more complex probability evaluation
see \cite{KnSZh91}.
Most of the examples given below involve broken line regression, where
a changes occurs in the slope of (the expected value of) a process, 
but the process itself
is continuous.  An important distinction between changes in level and
changes in slope is that the information about a change in level are 
usually found
in a relatively small neighborhood of the change-point, whereas information
about changes in slope can persist over long intervals.
A consequence is that confidence
regions for slope changes in broken line regression can be very large,
even though the existence of a slope change seems obvious.

Although the method of \cite{KnSZh91} produces a conservative result
under our assumptions, at least when $\rho = 0,$ here
we apply an asymptotic heuristic that is conceptually and 
technically simpler, and
can be adapted (Section 3.1) to deal with multiple signals.
The Kac-Slepian model process, for a mean 0
unit variance, twice mean-square differentiable, Gaussian 
process, say $U_s$,  suggests a
parabolic approximation for
the process in a neighborhood of $s_0$. (This result 
is very similar to the
local expansion of a quadratic mean differentiable log 
likelihood function, e.g.,
van der Vaart (1998)).
By a tedious calculation of means and covariances, we see that
\[ \Expec(U_s | U_{s_0}, \dot{U}_{s_0} ) =
U_{s_0} + (s-s_0) \dot{U}_{s_0} - \frac{(s-s_0)^2}{2} U_{s_0} \Expec_0[(
\dot{U}_{s_0})^2] + o_p((s-s_0)^2) \]
and 
\[ \Var(U_s | U_{s_0}, \dot{U}_{s_0} ) = O((s-s_0)^4)\].
These results and sufficiency (Proposition 5)  suggest that  
in a neighborhood of 
a putative signal $t$, where $Z_t$ is large we consider the
local approximation
\[
Z_{t+\delta} \approx Z_{t} + \delta \dot{Z}_{t} - \frac{\delta^2}{2} Z_{t} \Expec_0[(
\dot{Z}_t^2]).
\]
Maximizing over $\delta$,
we obtain
\begin{equation} \label{Kac-Slepian} \max_\delta(Z_{t+\delta}^2 -  Z_{t}^2) \approx \dot{Z}_{t}^2/ \Expec_0[(\dot{Z}_{t})^2], 
\end{equation} which is 
is approximately $\chi^2$ with 
one degree of freedom.   By inverting this relation, we
obtain as  an approximate  $1-\alpha$ 
confidence region the set of all $t$ that 
satisfy $Z_t^2 \geq \max_\delta Z_{t + \delta}^2 - \chi^2_{1-\alpha}$, 
where $\chi^2_{1-\alpha}$ is the $1-\alpha$ quantile of the $\chi^2$ distribution with one 
degree of freedom.  
This is in effect the same result one would obtain in large samples by 
inverting the log likelihood ratio 
statistic of a simple hypothesis if $Z^2$ were the 
log likelihood of a parameter satisfying standard regularity conditions.   
We conjecture that this approximation can be proved under 
additional regularity condtions, but here we use a small set of
simulations to suggest it is reasonable. 

A joint confidence region for $\xi$ and $t$ can also be obtained.  To the condition that
$Z_t^2$ must be within a given distance  of $\max_s Z_s^2$, we also require that 
$|Z_t - \xi \sigma(t)|$ is sufficiently small.  The pairs $\xi,t$ that satisfy 
$\max_\delta [Z_{t+\delta}^2 - Z_t^2] + (Z_t -\xi \sigma(t))^2 \leq c^2 $
provide a joint confidence region having the confidence coefficient

\begin{equation}
\Prob_{t, \xi}\{ \max_\delta [Z_{t+\delta}^2 - Z_t^2] + (Z_t -\xi \sigma(t))^2 \leq c^2\}
\end{equation}
\begin{equation} \label{jointconfidence}
= \Expec_{t,\xi}[ \Prob_{t,\xi} \{ \max_\delta [Z_{t+\delta}^2 - Z_t^2] \leq  c^2 - (Z_t-\xi\sigma(t))^2 
| Z_t\};
|Z_t -\xi \sigma(t)| \leq c]
\end{equation}
which equals the distribution of 
a  $\chi^2$ random variable with two degrees of freedom, again as if 
standard regularity conditions had been satisfied.  

In view of the heuristic elements of the preceding calculations, we give here
the results of a small simulation that suggests our approach is reasonable.
Consider a broken line regression model where observations have unit
variance and the only nuisance parameters are unknown intercepts and slope.  
The sample size is 100.
The slope is initally 0 and changes to $\xi$ at $t$.

\begin{table}[htp]
\caption{\label{simulations1} Confidence Sets for $t$.
Number of repetitions is $N= 1000$.}
\begin{center}
\begin{tabular}{|c|c|c|c|c}
\hline
$t$ & $\xi$ & Nominal Conf.& Empirical Conf.& Expected Size \\ \hline
30& 0.07 & 0.95 & 0.97 & 19 \\ 
50 & 0.05 & 0.95 & 0.96 & 30 \\ 
70 & 0.07 & 0.95 & 0.96 & 31 \\
30 & 0.05 & 0.90 & 0.89 & 21 \\
45 & 0.04 & 0.90 & 0.91 & 30 \\
60 & 0.03 & 0.90 & 0.88 & 49 \\
50 & 0.03 & 0.90 & 0.89 & 42 \\ \hline
\end{tabular}
\end{center}
\end{table}



\subsection{Numerical Examples}

\noindent
{\sl Example 1.  Extreme Precipitation in the United States.}
Extreme  precipitation in the United States is reported by the
National Oceanographic and Atmospheric Administration (NOAA) at
https://www.ncdc.noaa.gov/temp-and-precip/uspa/wet-dry/0.

This web site gives the area of the country where monthly rainfall 
exceeeded the  90th percentile 
of normal 
for the 125 years beginning in 1895.

Our test for at most one changes suggests a slope increase at 707 months,
in 1954,  
with a p-value of 0.002.
The value of $\rho$ is estimated to be 0.08.
The test applied to extremely low preciptation (area of the country below the 10th percentile)
is consistent with the hypothesis of no change in slope.  

\noindent
{\sl Example 2. Central England Annual Average Temperature.}
The average temperature in central England
from 1659 (\cite{Man74}) until 2019 is reported 
https://www.metoffice.gov.uk/hadobs/hadcet/cetml1659on.dat.  

\smallskip\noindent
The maximum $Z$ value is 
about 4.60, with $\rho$ estimated to be 0.16, and occurs
in 1969.  The appropriate (two-sided) p-value is about
$6\times 10^{-5}$.  
See Figure 1, where it appears that there might
also be earlier change-points.  
Although the mercury thermometer was relatively young in 1659,
those early years are interesting since
they appear to involve the ``little ice age.''  We return to these data below
for our discussion of segmentation.  Meanwhile it may be interesting to note that
an approximate 95\% confidence region based on the detected change in 1969 is (1895,1990),
which reflects the local shape of the plot in Figure 1 and suggests the possibility that
an increase in temperature may have actually begun during the late 19th century.

\begin{figure}
\begin{center}
\includegraphics[width = 5.0 in]{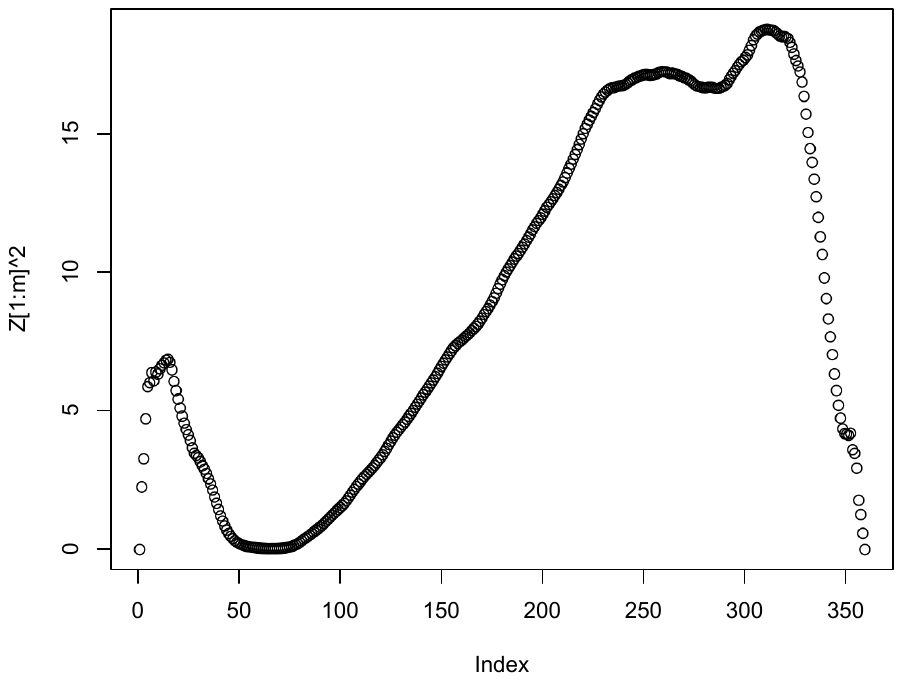} \\
\end{center}
\caption{Plot of $Z_t^2$ for central england temperature  data, 1659-2017.
\label{Fig 2}}
\end{figure}

Although we did not include a confidence region for Example 1, in that
case also, the plot of $Z_t^2$ exhibits two local maxima and a confidence
region substantially longer to the left than to the right of the 
maximizing value.

\subsection{Estimation of $\sigma^2$ and $\rho$: Simulations.}
Although we have assumed the variance known in our theoretical calculations,
in our numerical studies we have
estimated it by the residual mean
square under the null hypothesis.  Other possibilities 
when the data are independent, which may be less biased, are
sums of squares of first or second order differences:  $\sum_u (Y_{u+1} - Y_u)^2/2m$ or
$\sum_u (Y_{u+1} - 2 Y_u + Y_{u-1})^2/6m$.,
First order differences remove most of the effect of changing
mean values, and second order
differences also mitigate the effect of slope changes.

Regarding correlation of the observations, 
a useful model is one that helps to  control false positive 
errors with minimal loss of 
power, but otherwise estimating $\rho$ efficiently is unimportant.
Using a value of $\rho$ that is much too small can lead to an
increase in false positives.
But in the presence of change-points,
the maximum likelihood
estimator of $\rho$ under the null hypothesis 
that $\xi = 0$ can be upwardly biased, 
leading to a loss of power.

In Table 2 we consider the control of false positive errors 
when true and assumed first order autocorrelations differ.
The column labeled $\hat{\rho}$ contains assumed values, not
estimates. 

\begin{table}[htp]
\caption{\label{simulations2} Simulated False Positive Errors: 
$m = 100$, 0.05
threshold is  $b = 2.84$.  Number of repetitions is $N = 900$.}
\begin{center}
\begin{tabular}{|c|c|c}
\hline
$\rho$ & $\hat{\rho}$ & $\Prob\{\max_s |Z_s|\geq b\}$\\ \hline
0 & 0 & 0.033 \\
0.3 & 0 & 0.12 \\
0.3 & 0.1 & 0.076 \\
0.3 & 0.2 & 0.050 \\
0.5 & 0.3 & 0.043 \\
0.6 & 0.3 & 0.051 \\
0.7 & 0.3 & 0.059 \\
0.8 & 0.3 & 0.081 \\
0.9 & 0.35 & 0.049 \\
0.8 & 0.4 & 0.023 \\
0.9 & 0.4 & 0.036 \\ \hline
\end{tabular}
\end{center}
\end{table} 

The simulated example in Table 3 illustrates the problem
of estimating $\rho$ and the effect on the power to detect a change.  

\begin{table}[htp]
\caption{\label{simulations3} Simulated examples: $m = 150$, 0.05
threshold is  $b = 2.83$.
Locations of
break-points and changes in the slope occurring at the break-points are as indicated.
Rows where $\hat{\rho}$  has been set equal to $ \rho$ are included to 
illustrate the loss of signal 
that occurs because of bias in estimation of $\rho$.}
\begin{center}
\begin{tabular}{|c|c|c|c|c|c}
\hline
 $t$ & $\xi$ & $\rho$ & $\hat{\rho}$ & $\hat{\sigma}$ & $\max_s |Z_s|$\\ \hline
 --- & 0.00  & 0.50 & 0.48 & 0.98 & 1.60\\
 --- & 0.00  & 0.60 & 0.56 & 0.98 & 1.05 \\
 100 & 0.05 & 0.0 & 0.30 & 1.04 & 4.04 \\
 100 & 0.05 & 0.0 & 0.00 & 1.09 & 5.51 \\
 100 & 0.05 & 0.4 & 0.60 & 1.19 & 3.86 \\
 100 & 0.05 & 0.6 & 0.75 & 1.00 & 3.34 \\
 100 & 0.05 & 0.6 & 0.60  & 1.02 & 5.06 \\
  50 & 0.04 & 0.5 & 0.54 & 0.97 & 3.37 \\
  50 & 0.04 & 0.5 & 0.50 & 0.97 & 3.66 \\ 
  50 & -0.033 & 0.6 & 0.76 & 1.15 & 2.74 \\
  50 & -0.033 & 0.6 & 0.60 & 1.15 & 4.34 \\ \hline 
\end{tabular}
\end{center}
\end{table}

For the data in the third and fourth rows of Table 2, 
According to Proposition 2,
the expectation of $Z_t$ is 5.55, so the simulated value 
when it is known that $\rho = 0$ and $\sigma$ is estimated
to be 1.09 is very accurate.
In that case, if we 
use second order differences to estimate $\sigma^2$, the estimated value
of $\sigma$ is only 0.89 and the maximum $Z$ value is 6.77.   
In both rows 3 and 4, the  
maximizing value of $t$ is 94. 

It is apparent from Table 2 and other examples not shown that the 
existence of a
change-point causes the estimator of $\rho$ to be upwardly biased and results
in a loss of power compared to using the true value.  This
is usually worse if there are multiple change-points.  
If a plot of the data indicates long stretches without change-points, 
one might
use an estimator of $\rho$ from that segment.  
For the data in the next to last row of Table 1,
an estimate of $\rho$ based on the observations 65-150 yields 
the value 0.61, which
in turn leads to a maximum $Z$-value of 4.31. 
Although the theory developed here
does not justify this approach, numerous simulations 
suggest that
it works reasonably well. (See also the discussion of the Atlantic 
to use an arbitrarily chosen value of $\rho$ to pre-whiten the
observations and repeat the  analysis if that value appears to be
poorly chosen.
  
A second consideration in estimation of $\rho$ is robustness of our procedure 
against more complex forms of dependence.  Without exploring this question in
detail, consider again the next to last row of Table 1, but with 
a second order autoregressive coefficient of - 0.2 and $\xi = -0.06$.
A simulation of this case gives an estimator of 
$\rho$ equal to 0.6 based on all the data, with a  maximuom $Z$-value of 3.32,
and an estimator equal to 0.42 based on the second half of the
data, with a maximum $Z$-value of 4.85.  Other simulations, not given here,
suggest that at least for higher order autoregressive dependence, the first
order model we have suggested  provides reasonable protection 
against an inflated
false positive error rate.

\section{Segmentation.}
We now turn to the problem of segmentation when there may be several
change-points.  It is helpful conceptually to think of two different
cases.  In the first the signals exhibit no particular pattern.  In the
second two or more changes are expected to produce a 
definite signature.  For the problem of jump changes in levels of the process,
paired changes frequently take the form of an abrupt 
departure from followed by
a return to a baseline level.
Here we concentrate on the case of signals 
having no apparent pattern.  

The simplest method is binary segmentation: first we test for at
most one change, then iteratively partition the data at a putative change-point to test 
for additional changes
to the left or the right of the detected change-point.   A definite disadvantages of this method is 
that positive and negative changes in the same sequence of observations may cancel and lead
to a loss of power to detect any change.  In addition, for slope changes the score statistic has 
substantial local correlation, so it can occur that a second change is detected very close to the
first unless the interval searched is restricted so that its end-points are some
arbitrary distance away from already detected change-points.  

Next we 
consider versions of the two methods suggested in \cite{FaLSi20}
for detection of jump changes.
First we consider a pseudo-sequential procedure, called Seq below, 
where it is
convenient to write $Z(t,T)$ to denote the statistic $Z(t)$ when the interval
of observation is $[0,T]$.  We also use
a minimum sample size $m_0$ and
a minimum lag $n_0$, both of which 
we often set equal to 5.  Let $t_0 = m_0$ and let $T$ 
be the smallest {\sl integer} exceeding $t_0 + n_0$
such that for some $t_0  < t < T-n_0$, the value of
$Z(t,T)$ exceeds a threshold $b$ to be determined by a 
constraint on the probability of false positive error.  
The first change-point
is taken to be  $t_1$: equal to $ \arg \max |Z_t|$, 
or the smallest or largest value of $t$ for that smallest $T$. 
(In numerical experiments we have been unable to find a consistent preference.)
The process is
then iterated starting from $t_1$.  The nuisance parameter $\rho$ is
assumed to be known.  In practice it is estimated from the data or more
often from a
subset of the data to mitigate the effect of the bias
discussed above.  The other nuisance parameters, are estimated
from the current, evolving interval of observations, so that the nuisance 
parameters associated
with one change-point do not confound detection of another.  It would also be possible 
to estimate $\rho$ from the currently studied subset of the data, but this estimator appears
to be unstable.

To control the global false positive error rate, we want an approximation
for 
\[
Q = \Prob\{\max_{m_0 \leq t < T-n_0, T \leq m} |Z(t,T)| \geq b \},
\]
where $m$ is the number
of observations and $m_0,n_0$ represent minimal sample sizes, both of 
which we frequently take to be 5.  To emphasize that $T$ is a variable quantity,
we also write $\sigma^2(t,T),$  $\Psi(t,T)$, etc. 

\noindent
{\bf Theorem 1.} 
Let $\beta(t,T) = \{\Expec[V \partial V/\partial T)-\frac{1}{2} \partial \sigma^2(t,T)/\partial T\} / \sigma^2(t,T)$ and 
$\lambda_t = \Expec[(\dot{Z}_{t,T})^2].$  Assume that $\Delta$ and $\delta$ are positive
constants converging to 0.  Assume $b \rightarrow \infty$
and $b^2 \Delta $ converges to a positive constant, while $b \delta \rightarrow 0.$ 
Then
\begin{equation} \label{segmentprob}
Q \sim (2/\pi)^{1/2} b^2 \varphi(b) \int_{m_0}^{m}\int_{m_0 < t \leq T - n_0} 
\{\lambda_t^{1/2} \beta(t,T) 
\nu[b(2 \beta(t,T)\Delta)^{1/2}]\} dt dT.
\end{equation}

For a proof we 
also assume that $m$ is not too large relative to $b$, e.g.,
that $m = o(\exp(b^2/8))$, which ensures that asymptotically 
$Q$ is a small tail probability.  (See Remark (iii) below.)

For $t= i \delta$ and $ T = k \Delta$, we write 
the probability of interest as
a sum over $T$ of the sum over $t \leq T - n_0$  of 
the integral over $x>0$  
of the product of three factors:  (i) $\Prob\{ Z(t,T) \in b + dx/b\} \sim \varphi(b) \exp(-x) dx/b;$
(ii) $\Prob \{ \max_{j \geq 1} b[Z(t, T-j\Delta) - Z(t,T)] < -x | Z(t,T) = b+x/b \}.$ and
(iii) the conditional probability
that $[Z(t \pm \delta,T) - Z(t,T)] < 0$.  
The argument of (for example) Theorem 2.1 of \cite{FaLSi20}
shows that in the limit the conditional probability (ii)
converges to a probability of the form $\Prob\{ \max_{j \geq 1}
S_j < -x\}$, where $S_j$ is a Gaussian random walk with mean
equal in the limit to $\beta b^2 \Delta$ and variance 
$2 \beta b^2 \Delta$. The integral over $(0, \infty)$ 
of $\exp(-x)$ times this probability 
equals $\beta b^2 \Delta \nu [b(2 \beta \Delta)^{1/2}]$,
where $\nu$ is the function defined in Siegmund (1985) and
easily approximated numerically. 
By the Kac-Slepian approximation, for small
$s$, $Z(t+s,T) - Z(t,T) = s \dot{Z}(t,T) - s^2 \lambda Z(t,T)/2 
+ o_p(s^2)$, where $\dot{Z}$ is normal with mean 0 and variance
$\lambda$.  Hence
the conditional probability
in (iii) equals $\Phi(\lambda^{1/2} b \delta/2)]
- \Phi(-\lambda^{1/2} b \delta/2)$.  Since 
$b \delta \rightarrow 0,$ this last expression is asymptotic
to $\lambda^{1/2} b \delta/(2 \pi)^{1/2}.$  Putting these pieces together
and replacing the sums by Riemann integrals yields
(\ref{segmentprob}).  

\smallskip\noindent
{\bf Remarks and Examples.}
\noindent
(i) In (\ref{segmentprob}),  we 
treated the process $Z(t,T)$ as continuous in $t$
and have used $\delta$ and $\Delta$ as a differentials to 
replace sums by integrals..
In applications  
we always take $\Delta = 1,$
so one could leave the integral over $T$ in the form of a sum.
Numerically this makes essentially no difference.
We could also take $\Delta$ to
0 faster than $1/b^2$, which would effectively treat the 
process in $T$ as a continuous, but not differentible, 
process.  This leads
to a very conservative approximation if the process is
approximately Gaussian, but presumably more robust
to departures from normality.
(ii) For a numerical example, suppose $m = 365$, $b = 4.0 (4.41) $,
$m_0 = n_0 = 5.$  Then (\ref{segmentprob}) gives the 
approximation 0.05 (0.01).  When $T$ is allowed to vary 
continuously, the threshold corresponding to  0.05 (0.01) 
is 4.43(4.80.)  A small simulation with $m_0 = n_0 =5$ and a sample
size of 10000 to test the
accuracy of the first set of thresholds gave the values 0.0493 
and 0.0111.

\noindent
(ii)  It is possible to make the preceding 
argument rigorous, but the many details that one can see 
in related arguments
suggest that a complete proof  would be onerous and not add any new 
insights (cf.  \cite{SY2000}, 
\cite{FaLSi20}).  

\noindent
(iii) The  condition that $m$ not be too large compared to
$b$ guarantees that there is asymptotically at most one excursion
above a high threshold $b$, so the right hand side (RHS) of  
(\ref{segmentprob}) converges to 0, and hence our 
approximation lies
in the domain of small tail probabilities.  For much 
larger $m$, there can be several excursions above $b$, 
so RHS(\ref{segmentprob}) may be bounded away from 0,
and one expects in that case that a suitable  approximation
has the form $Q - [1- \exp[-RHS(\ref{segmentprob})]
\rightarrow 0.$  See \cite{SY2000} for the global arguments
for a similar example. 

\noindent
(iv) When the process $Z(t,T)$ exceeds a suitable
threhold, there is typically a cluster of values of $T$
and $t$ 
where $|Z(t,T)|$ exceeds  
that threshold.  In practice 
we might for the smallest such $T$  choose as an 
estimate of the position of a slope change 
the smallest 
$t$, the value 
that maximizes $|Z(t,T)|$, or perhaps some other
value. One hopes that experience 
eventually suggests an optimal
choice of $t$, but that does not seem to be the case.   
Examples where
there is only a small number  of values of $T$ and $t$
with the values of $t$ only slightly smaller than $T$
may indicate the presence
of outliers and suggests that we  skip to a  
larger value of $T$.
Once a value of $t$ is selected, we then iterate the
process starting from that $t$, and 
the background linear regression 
is estimated using the data beginning with that $t$.  

\noindent
(v) We have  used the term ``pseudo-sequential'' 
since most applications discussed in this
paper involve a fixed amount of data.
However, we estimate $\alpha$, $\beta$ and $\sigma^2$  using
the data since the last change detected, so the 
method can also be applied sequentially.
See Section 4 for an application to detecting upswings in COVID-19..  

\medskip
A second method of segmentation, called MS below 
(for Maximum Score Statistic),  uses (in obvious notation) 
\[ 
\max_{m_0 \leq T_0 < t < T_1 < m} Z_{T_0, t, T_1}.
\]
An asymptotic approximation to the null probability that this expression exceeds $b$
can be computed as a corollary of the  
argument behind Theorem 1, and we obtain the approximating expression
\[
(\frac{2}{\pi})^{1/2} b^4 \varphi(b) \times
\]
\begin{equation} \label{segmentprob2} \sum_{m_0\leq T_0 < T_1 \leq m}\int_{T_0 <  t < T_1}
(\lambda_t)^{1/2} \beta_0(t,T_0)\beta(t,T_1)
\nu\{b[2 \beta_0(t,T_0)]^{1/2}\}\nu\{b[2 \beta(t,T_1)]^{1/2}\} dt,
\end{equation}
where $\beta_0(t, T_0)$ is defined similarly to $\beta(t,T)$.  The evaluation 
of (\ref{segmentprob2})
can be simplified by a summation by parts and the observation that the various functions of 
three variables, $(T_0, t, T_1)$, occurring in (\ref{segmentprob2}) are
in fact functions of two variables:  $(t-T_0, T_1-T_0)$.  

This method determines a list of putative local signals together with
generally overlapping  backgrounds.  We
consider different algorithms for selecting a set of local signals 
and backgrounds from this
list, with a preference for the additional  constraint that no 
background overlap two local signals.
In \cite{FaLSi20} examples of particular interest were the shortest background, 
(cf. also \cite{BCFr16}), and the largest $Z$ value.
Since we have paid in advance
for false positive errors,  we can also consider 
other methods for searching the list, e.g., subjective 
consistency with a plot of the data.  We also find it 
convenient computationally 
to follow the suggestion of \cite{Fr14} by searching a 
random subset of intervals, which seems to work very well.
When it appears that the number of changes is small, 
as in the examples considered in this paper,
one might also use a combination of methods, e.g., searching a small
random subset of intervals at a low threshold to generate candidates,
then sorting those by hand with tests to detect at least one change
at the appropriate higher threshold. 

The methods Seq and MS  ``bottom up'' methods in the sense that they attempt to
identify one local signal at a time against a background appropriate for that signal--an 
advantage if there is heterozygosity. 

A popular ``top down'' method to detect abrupt level changes scans all the data looking for
a pair of change-points
(e.g., \cite{OlVe04}).  In some cases, especially if there is only one change to be
detected, the method will ``detect'' two changes, as required, but it will put one of them near
an end of the data sequence, where it can be recognized and ignored.    
After an initial discovery 
of one or two change-points, the sequence is
broken into two or three parts by those change-points, and the process is iterated as long
as new change-points are identified. 
Consider the statistic  $\max_{s < t-h} R(s,t)$,
where
\begin{equation} \label{2locusstatist}
R^2(s,t) =  (Z_s, Z_t) \tilde{\Sigma}_{s,t}^{-1} (Z_s,Z_t)',
\end{equation}
$\tilde{\Sigma}$ is the covariance matrix of $(Z_s, Z_t)'$,
and $h$ is a parameter that represents a minimum distance between changes that
we find interesting (usually taken to be 5 or 10 in the examples below).
An appropriate threshold may  be determined
from an approximation similar to (9) although the details are
substantially more complicated, since $R(s,t)$ is related to, but is not
a Gaussian process.  

Let $\rho = \rho(s,t)$ denote the correlation of $Z_s$ and $Z_t$ and put
$U(\theta,s,t) = \cos(\theta) Z_s + \sin(\theta)(Z_t - \rho Z_s)/(1-\rho^2)^{1/2}$. 
Then $U$ is a
Gaussian process and $R(s,t) = \max_\theta [U(\theta,s,t)].$
Now a formula similar to (9) is applicable.  Computation of the final result 
is complicated by the fact  that the integral is three
dimensional and involves the determinant of the $3 \times 3$
covariance matrix of  $\dot{U} = (U_\theta, U_s, U_t)'$, where
subscripts denote partial derivatives. The result is

\noindent
{\bf Theorem 2.} As $b \rightarrow \infty$,  
\begin{equation} \label{2locusapprox}
\Prob\{ \max_ {s<t - h} R(s,t) > b \}  \sim [b^2 \varphi(b)/(2\pi)^{3/2}] \int_{0 <\theta < 2 \pi, s < t-h} {\rm det}[\Expec(\dot{U} \dot{U}')]^{1/2} d\theta ds dt.
\end{equation}


\subsection{Confidence Regions Revisited.}
We can also use 
(\ref{2locusstatist}) to find confidence regions for two (or more) change-points
$t_1 < t_2$.
First note that with $t = (t_1, t_2)$, and $\xi = (\xi_1, \xi_2)$, a
straightforward calculation extending the result given in Proposition 2
shows that
$\Expec_{t,\xi}(V_s) = \sum_i [\xi_i \sigma(t_i) \Sigma(s,t_i)]$, where $\Sigma(s,v)
= \Cov(V_s, V_v).$  This produces by  a second
calculation the
conclusion that the
log likelihood of $V_s, s = 1, \cdots, T$
equals $\sum_i \xi_i \sigma(t_i) V_{t_i} - .5 \sum_{i,j} \xi_i \xi_j \sigma(t_i) \sigma(t_j)  \Sigma(t_i, t_j)$.
Hence the log likelihood ratio statistic, given $t$, for testing that
$\xi$ is the 0 vector equals
\begin{equation} \label{chisqstat}
 (V_{t_1}, V_{t_2}) \Sigma(t_1,t_2)^{-1} (V_{t_1}, V_{t_2})'
= || \tilde{Z}_{t_1,t_2} ||^2,
\end{equation}
where $\tilde{Z}_{t_1,t_2} = \Sigma(t_1,t_2)^{-1/2} (V_{t_1}, V_{t_2})'$
has a standard bivariate normal distribution when $\xi = 0$, and
$|| \cdot||$ denotes the Euclidean norm.  It follows from
sufficiency and the Kac-Slepian argument used above that
\[ \max_{\delta_1, \delta_2} [|| \tilde{Z}_{t_1+ \delta_1,t_2 + \delta_2} ||^2
-|| \tilde{Z}_{t_1,t_2} ||^2]
\]
has a $\chi^2$ distribution with two degrees of freedom and can be used as
above to obtain a joint confidence region for $(t_1, t_2)$.

A joint confidence region for $t, \xi$ follows by an argment similar to that given above.

\subsection{Applications.}  In this section we consider a number of 
data sets, where
broken line regression may be an appropriate model.  
Although it may not fit the data globally as well as other possibilities, it
suggests questions regarding the times of sharp changes of direction
and reltively small
changes in long term trends.  Since all these examples  involve
observational time series, an assumption of independence seems 
inappropriate, although in several cases the estimated 
value of $\rho$ is small enough
to be assumed equal to zero. 

As a step in our evaluation of  the segmentation methods of Section 3,
we assume that the
detected changes are in fact correct,so  
our model becomes a standard linear model, with $\rho$ 
the coefficient of a one-step lagged regressor.  Estimation of 
the regression parameters of that model with the change-points assumed 
known provides an idea of the size and importance of different detected changes,
the adequacy  of the model as judged by the value of $R^2$ 
and $\hat{\rho}$, and 
the reasonableness of our
choice for the value of $\rho$ in our initial segmentation.  In cases 
where the number
of changes may be in doubt, e.g., when the method MS detects fewer changes, 
probably 
because of its stricter significance standard, we also compute a BIC score 
for the
model with the larger number of slope changes compared to the model 
with fewer slope
changes.  Although the broken line model is an irregular statistical model, according to
calculations in the unpublished Ph. D. thesis of Yi Liu, 
the customary assessment of 
a penalty by counting parameters (multplied by the log sample size) is the same as for a 
regular model. 

In the figures that follow, we have plotted the
data superimposed on the results of a multiple
regression analysis.  For ease of viewing
the autocorrelation is
not incorporated into the plot, although it is included in the
evaluation of $R^2$.  

\medskip\noindent
{\sl Example 3: World Average Land and Ocean Temperature 
Anomalies: 1850-2020}

Temperature anomaly time series are calculated by the Hadley 
Center in
Great Britain, NOAA and the Goddard Institute for Space Sciences (GISS)
in the United States, and by Berkeleyearth.
Presumably approximately the same data from the same instruments are available
to all groups,
which they then process by different methods.  
The main obvious difference
is the beginning date for the different time series: 1850 for the 
Hadley Center, 1880 for NOAA and GISS, and varying times for different countries, 
some as early as about 1750 for Berkeleyearth. 

The data we consider in this example are the Hadley Center 
time series for land and ocean temperature  
anomalies, which begins in 1850 and contains
171 observations.  
Using MS with $\rho = 0.2$, we detected three slope changes,
an increase in approximately 1910, a decrease in 1942, 
and a final increase in 1977.   
The output of a linear
analysis using those three changes is displayed in Figure 2.
Since the data appear heteroscedastic, with more variability at the beginning of the
series, in our linear model we used a weighted analysis, with weights derived from the 
uncertainty measures given for the data.  
The value of $R^2$ is estimated to be 0.89, and the value of
$\rho$ is estimated to be about 0.3.  With this value of $\rho$, both Seq and (\ref{2locusstatist}) suggest 
essentially the same three changes.   MS and Seq, which detect changes against a local background, 
give roughly the same estimate of $\sigma^2$, which suggests that 
heterozygosity is not a serious problem for these data.

We discuss some related series and the problem of heterozygosity  in the Online Supplement.

\begin{figure}
\vspace{-1in}
\begin{center}
\includegraphics[width = 4.0 in]{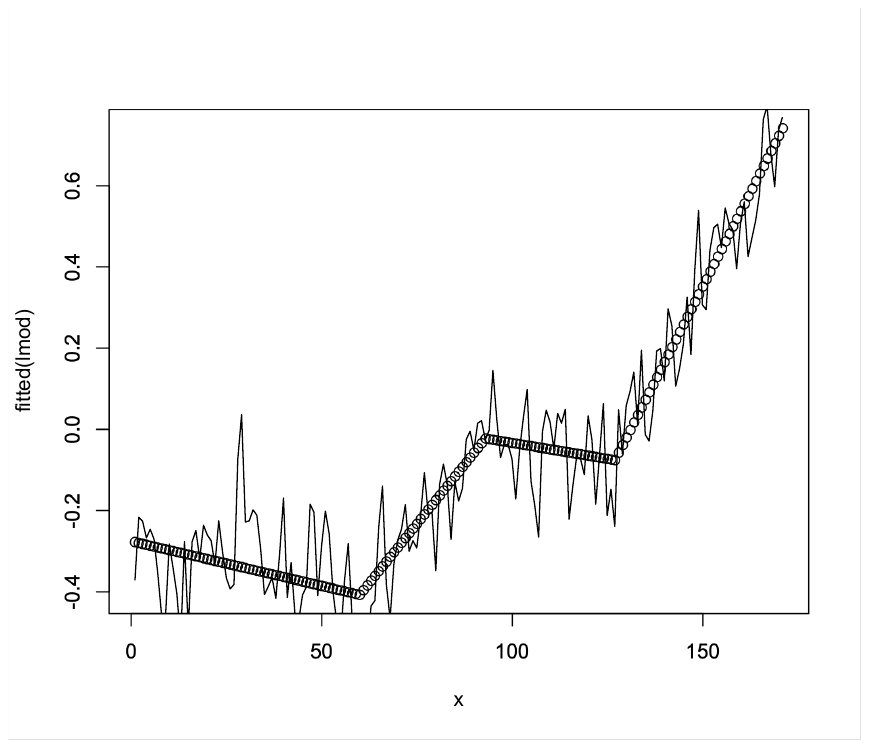} \\
\end{center}
\caption{Linear Analysis of Annual World Land and Ocean
Anomalies: 1850-2020. \label{Figure 2}}
\end{figure}

\medskip\noindent
{\sl Example 4. Excess Deaths in the United States.}
We consider next the number of excess deaths in the U.S beginning with 
the first week of
January 2020.  To avoid the problem of late reporting of deaths, we used 
86 weeks, although 93 were available when this manuscript was written. 
An interesting feature of excess deaths is that 
on the one hand they account both for deaths from COVID-19 and
deaths that occurred because appropriate medical care was unavailable due to the 
pandemic, and on the other they
do NOT count deaths from influenza that did not materialize 
because responses to the pandemic.

For these data   
(\ref{2locusstatist}) with $\rho = 0.3$   
detected eight slope changes.  Seq missed the first slope change, but otherwise agreed
to within one week on the other detections.
Results displayed in Figure 3 show a linear analysis 
with the changes suggested by (\ref{2locusstatist}) with some visually suggested adjustments
to agree with Seq.  For this linear analysis $R^2 \approx 0.97$ and
$\rho \approx 0.19$.  

\begin{figure}
\vspace{-1in}
\begin{center}
\includegraphics[width = 3.7 in]{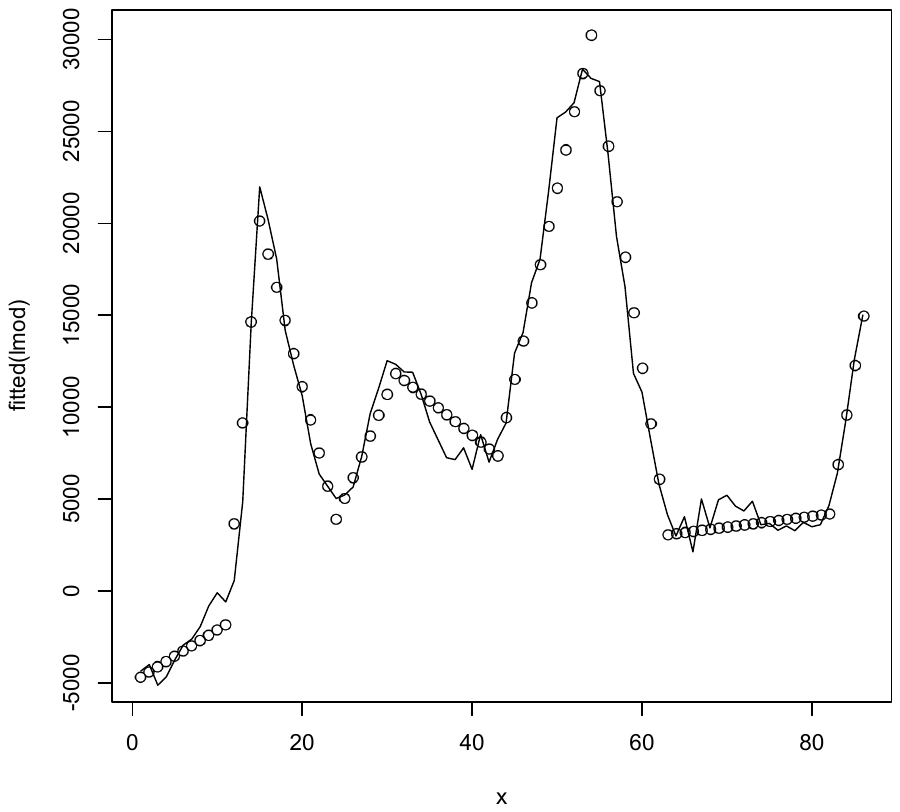} \\
\end{center}
\caption{Linear Analysis of Excess Deaths in U.S. for 86 weeks from 01/01/2020. 
\label{Figure 3}}
\end{figure}

\medskip\noindent
{\sl Example 2 revisited: Central England Average Annual Temperature.}
A challenging case for segmentation is the central England temperature data 
for 361 years beginning in 1659.  A plot of the data suggests 
possible changes very
early in the series with a much larger change 2-3 hundred 
years later. 
The first 80-100 years may be less reliable, due to the relatively primitive thermometers
available at that time, but the early data are nonetheless interesting in view of the
``little ice age,'' which overlapped these years.  As reported above, a test of 
at least one change produces a significant increase in slope in 1970 
with  an estimated  autocorrelation 
of 0.159, and a confidence interval going  back about a hundred years,
reflecting the increase in the plot of $Z_t^2$ in Figure 2 that begins  in the late
19th century.  If we use Seq with the threshold 4.4 and the autocorrelation 0.1, we 
obtain a broken line that decreases steeply for the first 33 years, then increases 
for about 16, bringing it back more or less to where it started.  It then has a slightly
negative slope until about 1890, when it begins to increase.  No further increase in the
20th century is detected, which is perhaps surprising given our results when testing for
at most one change and an accelerated rate of increase starting about 1970 for many other
temperature series.
For MS with a threshold of 4.81 and the same value for $\rho$, we
find evidence for a change in 1890 or in 1983, without a compelling reason to choose one of
the other.. 
Figure 3 suggests that either (or both) could be correct.  
A linear analysis indicates that
the model with early changes in 1692 and 1708, then a later change in 
either  1890 or 1983 (not both) leads to 
an $R^2$ of approxmately 0.31 and a value of $\rho \approx 0.1$.
The statistic (20) to detect two changes at a time detects  
changes in  1693, 1708, and 1978.  BIC has a prefernce for three slope changes over four
and a very small preference for the change in 1890 over the change in 1983.  

\medskip\noindent
{\sl Example 5:  Age specific suicide rates in the United States.}
A small sample size example is suicide rates in the US from 1990 through 2017,
which can be found on the web site
ourworldindata.org.  
Autocorrelation appears to
be very small, so we assume it is 0.   All methods detect a 
substantial slope increase about 2000.
For MS this is the only change detected, and it produces an $R^2$ of 0.91.
Seq detects in addition, a second increase about 2014.
The statistic (\ref{2locusstatist}) detects three  changes,
a decrease in slope about 1995 before a large increase in 1999, and
a third increase in 2014.
A linear analysis with all three changes yields an $R^2$ of 0.97.

\section{Sequential Detection.}
In this section we consider using Seq for online detection of slope changes,
which for independent normally distributed observations has been studied
using different metholds by \cite{KPY03}.
We estimate the nuisance parameters locally, so that 
when we have observed $Y_S, \ldots, Y_T$, we
use only those values to estimate the nuisance parameters $\alpha$, $\beta$, 
and $\sigma^2$. 
Since estimates of $\rho$ are unstable when sample sizes are small, here
we set $\rho = 0$ and consider the possibility of reanalyzing the
data with a different value
if autocorrelation appears to play an important role.

We begin by considering the data for patient B in 
Smith and Cook (1980), one of a
beautiful series of papers 
by Smith and
co-authors concerned with monitoring renal transplants for
rejection.  
The data are given graphically in the paper, and by eye we extracted the
10 daily values 35,45,49,64,75,71,69,60,31,21.  Using Seq 
with the 0.05 threshold $b = 2.56$, we find after 8 observations that
a change has been detected and estimate the change to occur with the
6th observation.  A linear analysis with this change-point indicates
that the autocorrelation of 0 is reasonable and $R^2 \approx 0.92$.

Several features of the COVID-19 pandemic provide 
interesting examples where broken line 
regression might be useful. 
First we consider the possibility of 
sequential detection of a a slope increase in daily incidence 
in the initial phase of the pandemic and
allow a probability of 0.05 (or 0.01) of incorrectly detecting a slope increase 
within one year.  
By (14) with $m_0 = m-m_1 = 5$, 
this produces a (one-sided)
threshold of $ b \approx 3.81 (4.24)$. 

Data are from {\sl ourworldindata.org}.

As a first example we consider South Korea, which experienced its first
case on 20 January, 2020.  On 21 February, both the 0.05 and 0.01 
thresholds were crossed and the change attributed to 6 days earlier.

A second example is Italy, which experienced
its first cases on 31 January.  
Using either the 0.05 or the 0.01 criterion, we detect a slope 
increase on 26 
February and estimate
the change to have
occurred about 5 days earlier.  
A ``lockdown'' of the entire country appears to have 
occurred on 11 March.  

\smallskip\noindent
{\sl Remarks.} 
A potential difficulty with the sequential application of our model in the initial
phase of an epidemic is that detection of an increase in slope may occur while the sample 
size is still small, so using a normal approximation may be anti-conservative.  An alternative 
leading heavier tail probabilities would be to assume a Poisson model
with a constant mean value that may at
some point start to increase linearly.  Different possibilities  for
development are  (i) to extend our theoretical 
calculations, 
(ii) simulate a threshold to check and perhaps adjust results based on an 
assumption of normality, or (iii) do a linear analysis of a log linear Poisson model 
to confirm the results.  For the two examples given above, informal analysis confirms the
results given.  In other examples, detection was delayed by one or two days.

The sequential method may also be helpful in detecting second and 
third waves of infections following a relatively quiescent interval.
We consider Santa Clara County, California, which
like localities in the United States
experienced a very large number of cases in early 2021, 
followed by a sharp
decline in daily incidence during May.
We use data from  {\sl data.scc.org.gov}.
On 1 June, after a month when the daily 
incidence had been decreasing to fewer than 100 cases a day
in a population of about two million, we started Seq with
$\rho = 0.5$, which is approximately the value given by a 
complete segmentation of 600 observations.  The result with a 
(two-sided) 0.05 threshold  was detection
of a slope increase in the middle of July that was estimated to have begun 
during 
the last week of June.  
Using the 0.01 threshold delayed detection by 6
days.  

\section{Discussion.}
We have studied score statistics 
to detect local signals
in the form of changes of level, slope, or (in the
online Supplement) the autoregressive coefficient.  
To  segment the observations, we consider two ``bottom up'' methods 
patterned 
after the methods 
of multiple change-point segmentation 
in \cite{FaLSi20} and
one new ``top down'' method, defined in (\ref{2locusstatist}).
The three methods have different strengths and weaknesses.  Although both 
bottom up methods can be applied algorithmically, it appears to us that
satisfactory segmentations often
require some judgment regarding the appropriate number of 
changes and their locations.  The method Seq can be applied 
sequentially if it is important to detect a change in slope quickly
after its occurrence.    

Asymptotic approximations for  false positive error probabilities for the
three methods appear to be new.  For the 
bottom up methods the approximations  involve random fields
that are differentiable in one coordinate, but not differentiable
in others.  Both the statistic that tests for at most one change 
and (\ref{2locusstatist}) can be used to obtain confidence regions.  

Estimation of the nuisance parameters, $\rho$ and $\sigma^2$, pose 
special problems.
For the bottom up methods, in examining an interval of observations for possible change-points,
intercept, level, slope and variance  of the process
are estimated locally, i.e., using only the data from the interval under consideration.
For the autoregressive coefficient, our theory suggests 
using the (global)
maximum likelihood estimator
estimator under the 
hypothesis of no change.  Since this estimator can be badly biased
and result in a loss of power when there are change-points, 
especially multiple change-points,
we also consider {\sl ad hoc} methods based on
estimation from different parts of the data or on an arbitrary choice 
that it may be appropriate to change after subsequent analysis.
 

The multiple regression analsis suggested in Section 3.1, based on the 
assumption that detected break-points
are correct, allows us to see if our segmentation is reasonable, estimate the magnitude of the
detected changes, and reconsider our chosen value of $\rho$,  
if it appears that our original choice was inappropriate.  

Since the value used to estimate $\rho$ seems to have little
effect on the accuracy of the location of detected change-points,
a different approach would be initially to use the value $\rho = 0$,
or some other small value, 
apply the multiple regression
analysis and then iterate the process with
the value of $\rho$ suggested by the regression analysis.  
We have not chosen this path
because our goal has been to use a value of $\rho$ obtained with 
at most a small amount of  ``data snooping,''  in order
to be comfortable in saying  that due to minimal re-use of the data the 
false positive error rate has been
adequately controlled. 

In our view the
most interesting 
potential applications of our methods
are those where there is evidence of a local signal that may ask for an explanation.
For example, the ubiquity of rapidly increasing tempertures since about 1980
coupled with increased heat trapping gases leaves little room to doubt global
warming.  But some examples suggest that global warming began in late 19th or early
20th century, although that evidence is neither as strong nor as consistent.

We also find examples where
our methods for detecting multiple slope changes lead
to a reasonable regression fit without suggesting that 
individual changes require explanations. 


For problems involving sparse, bump like changes as 
discussed briefly in the Supplement, some 
superficial analysis
suggests that our methods
work well, provided that (as we assumed in simulations)
most of the data are consistent with the null model, 
and hence global
estimation of nuisance parameters does not pose a serious problem.
In view of the ubiquity of genomic applications where ``bumps'' 
have different signatures suggested by a combination of science and experimental technique, 
and where the background may involve use of a control group it seems
worthwhile to pursue a more systematic study.

Here we have considered signals in one-dimensional processes, where the number of possible 
``shapes'' of the signals is relatively small.  In view of the much larger
variety of possible multi-dimensional signal shapes and the variety of approaches
already existing in the literature for these problems, a 
systematic comparative study
would be valuable.

\bibliographystyle{Chicago}

\clearpage
\begin{center}
{\large\bf SUPPLEMENTARY MATERIAL}
\end{center}

\bigskip
\centerline{\bf Appendix A:  Some Basic Calculations}
\medskip
In this appendix we record basic formulas we have used in the paper.  Assume the model is
given by (2) with $M = 1$ and $\mu_u = \alpha + \beta[(u-(T+1)/2)/T]$.  The expressions given
below are asymptotic for  $T>> 1$;
and for some specific examples, the results used are computed
in continuous time.  Hence, for example, $\sum_{u=1}^T f[(u-t)/T]$ may be evaluated as
$T \int_0^1 f(u-t/T)du$, and $\Expec(Y_u)$ may be computed after solving the differential
equation (3).  In the change-point problems of (\cite{FaLSi20}), where the
result in considerable loss of accuracy, but it seems much less consequential here where 
for many of our examples, in particular for broken line regression, the stochastic
processes under consideration have piecewise differentiable sample paths.
Let $\gamma = 1- \rho$.  Then 
\begin{equation}
Y_u = (\sum_{j \geq 0} \rho^j)[\alpha + \beta(u/T-1/2]) + \sum \rho^j \epsilon_{u-j} + O_p(1/T), 
\end{equation}
or alternatively the solution of (3) with $Y_0 = 0$ is
\begin{equation} \label{integralrepresentation}
Y_u =  \int_0^u \exp[-\gamma(u-s)] \{(\alpha + \beta (s/T-1/2) +\xi f[(s-t)/T])ds + dW(s)\}.
\end{equation}
Hence, under the null hypothesis that $\xi = 0$ we have
\[ \Psi(t)' = \Expec(\ell_\xi \ell_\theta)' \sim T(1,1/12, (1-\rho)^{-1}[\alpha  +
\beta /12].  
\]

Let $\Delta_u = Y_u - \rho Y_{u-1} - \alpha - \beta(u-T/2).$
The ingredients of the representation $Z_t  = [\ell_\xi(t) - \Psi(t)'A \ell_\theta]/\sigma(t)$  
introduced in
Section 2 are
\[  \ell_\xi(t)  = \sum_u \Delta_u f[(u-t)/T];
\]
\[ \ell_\theta = (\sum_u \Delta_u, \sum_u \Delta_u [(u-(T+1)/2)/T], \sum_u \Delta_u Y_{u-1})';
\]
\[ \Psi(t) \sim (\sum_u f[(u-t)/T], \sum_u f[(u-t)](u/T-1/2), \sum_u \{f[(u-t)/T] (\alpha+ \beta[(u-T)/(2T)]\}/(1-\rho))'.
\]
For the special case of broken line regression in continuous time the first and
second coordinates are respectively
$g_1(t,T) = (T-t)^2/2T^2, \; g_2(t,T) = [(T-t)^3/12 - t(T-t)^2/4]/T^3$.

The matrix $A^{-1} = \Expec(\ell_\theta \ell_\theta')$ is straightforward to compute and somewhat
tedious  to invert.
A very useful result is
\begin{equation} \label{equ}
\Psi(t)' A = ( \sum_u f[(u-t)/T], \sum_u f[(u-t)/T][(u-1/2)/T]/\sum_u(u/T-1/2)^2, 0).
\end{equation}
It follows that the numerator of $Z_t$, 
has covariance function 
\begin{equation} \label{covariance}
\Sigma(s,t)  \sim \sum_u f[(u-s)/T]f[(u-t)/T] - \Psi(s)' A \Psi(t),
\end{equation}  
which does not depend on $\alpha$, $\beta$, nor on $\rho$. 
In particular 
\begin{equation} \label{variance}
\sigma^2(t) \sim \sum_u f[(u-t)/T]^2 - \Psi(t)' A \Psi(t).
\end{equation}

To derive (\ref{segmentprob}) and (\ref{segmentprob2}) (see below), we must consider $T$ as variable and study
the quantities $\sigma^2, \Psi,$ etc., as functions of both  $t$ and $T$.  
It is obvious from (\ref{variance}) that $\partial \sigma^2(t,T)/\partial T$
does not depend on nuisance parameters.
Since $\Psi(t,T)$ depends on the nuisance parameters $(\alpha, \beta, \rho)$ only
in its third coordinate, if we first differentiate (\ref{equ}) with respect to $T$,
then multiply on the right by  $\Psi(t,T)$ it follows from (\ref{equ}) that neither 
$(\partial \Psi/\partial T)' A_T \Psi(t,T)$ nor 
$\Psi(t,T)' (\partial A_T /\partial T) \Psi(t,T)$ depends on nuisance paramters.
The approximation (\ref{segmentprob2}) requires that we introduce $T_0 < t$ as
a third parameter.  For broken line regression we also require
\begin{equation} \label{equ2}
\Psi' A \partial \Psi/\partial{T_0} = -3T^3(1-t/T)^2 g_2[(t-T_0)/T,(T_1-T_0)/T]/(T_1-T_0)^3.
\end{equation}

To study the behavior of $Z_t$ when $\xi$ is different from 0, and to
help with the interpretation of a plot of $Z_t$, we put
$t = (t_1, \ldots, t_k)$, $\xi = (\xi_1, \ldots, \xi_k)$ and use the notation
$\Expec_{t,\xi} (\cdot)$.  By combining the results given above we find that
\begin{equation}
\Expec_{t, \xi} (Z_s) = \sum_j \xi_j \Sigma(t_j,s)/\sigma(s).
\end{equation}

\centerline{\bf Appendix B: Other Examples.}
This appendix contains a number of miscellaneous examples.  We begin with several related to
climate change.

\medskip\noindent
{\sl Example 5.  Extreme Heat.}
Again we use the NOAA website, where extreme heat is defined by the area of the country where the
temperature is above the 90th percentile.  Data are for 1518 months beginning in 1895.  Using
Seq with $\rho = 0.1$, we detected a slope decrease in about 1928 and a larger increase in 1974.  
A linear analysis gives estimated values of $R^2 \approx 0.44$ and $\rho \approx 0.11.$
We used annual averages and with (\ref{2locusstatist}) also detected two changes, in 
the years 1934 and 1974.
For the annualized data, a linear  analysis prefers the changes 
suggested by  (\ref{2locusstatist}).

To provide a rough description of a two-dimensional  confidence region for 
these data, we can first fix the first 
change-point at 1934, and then a $90\%$ region would extend from 1964 to 1990. 
If we set the upper
change-point at 1975, the $90\%$ region would extend from 1911 to 1965.  
The longer interval around 
the lower change-point reflects the fact that the size of the slope 
change at 1934 is only about 
4/5 as large
as the size of the slope change at 1975.  See Figure 5.

\begin{figure}
\begin{center}
\vspace{-1.0in}
\includegraphics[width = 5.0 in]{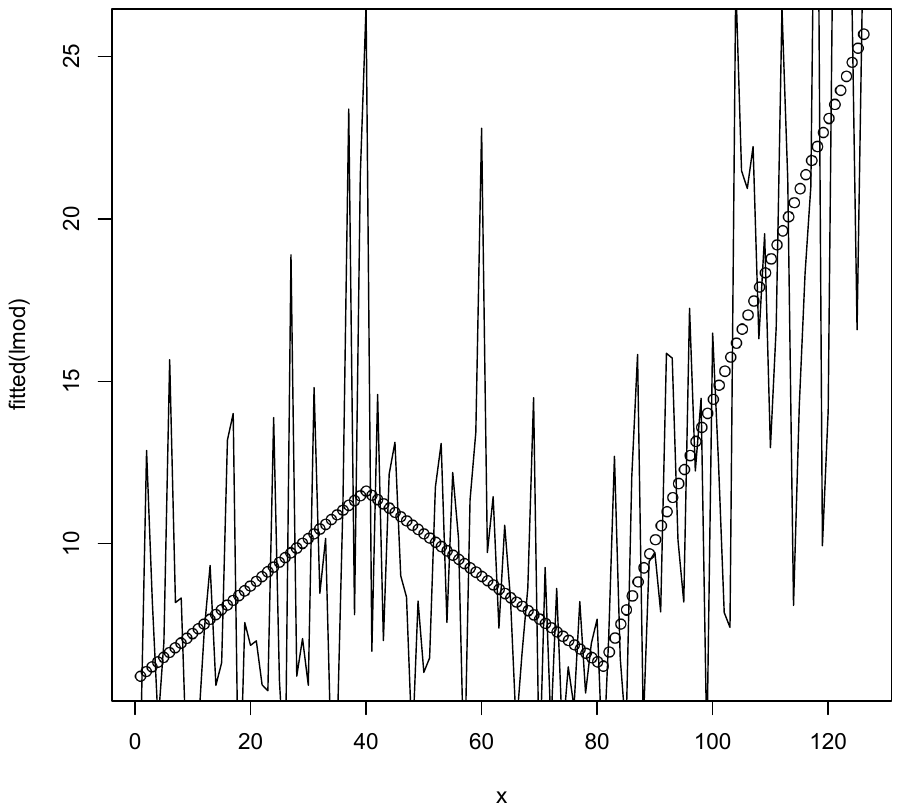} \\
\caption{Linear Analysis of Extreme Heat. \label{Figure4}}
\end{center}
\end{figure}

\medskip\noindent
{\sl Example 6: Ocean Temperature Anomalies.}
An interesting example is provided by time series of ocean temperature anomalies.
Both NOAA and the Hadley Center provide relevant data.  We consider here global annual 
anomalies from the Hadley Center, with 171 observations from 1850 and 2020.  Since a measure
of uncertainty is provided and early values appear to be less accurate than recent
values, we use a weighted analysis.

Using the value $\rho = 0.3$, both Seq and (\ref{2locusstatist}) agree on slope  changes 
in 1878, 1911, 1952, and 1972.  A linear analysis suggests that the first change is 
superfluous.  The output of a linear analysis, with the first change included, appears in Figure 5.

For MS we used an unweighted analysis, but selected slope changes on the basis of the shortest
background.  This resulted in obtaining the second, third, and fourth changes of the 
preceding analysis,
with the mid 20th century change put at 1944.  This is close to the 
outcome that we obtained
in Example 3 for land and water together.  

\begin{figure}
\begin{center}
\vspace{-1.0in}
\includegraphics[width = 5.0 in]{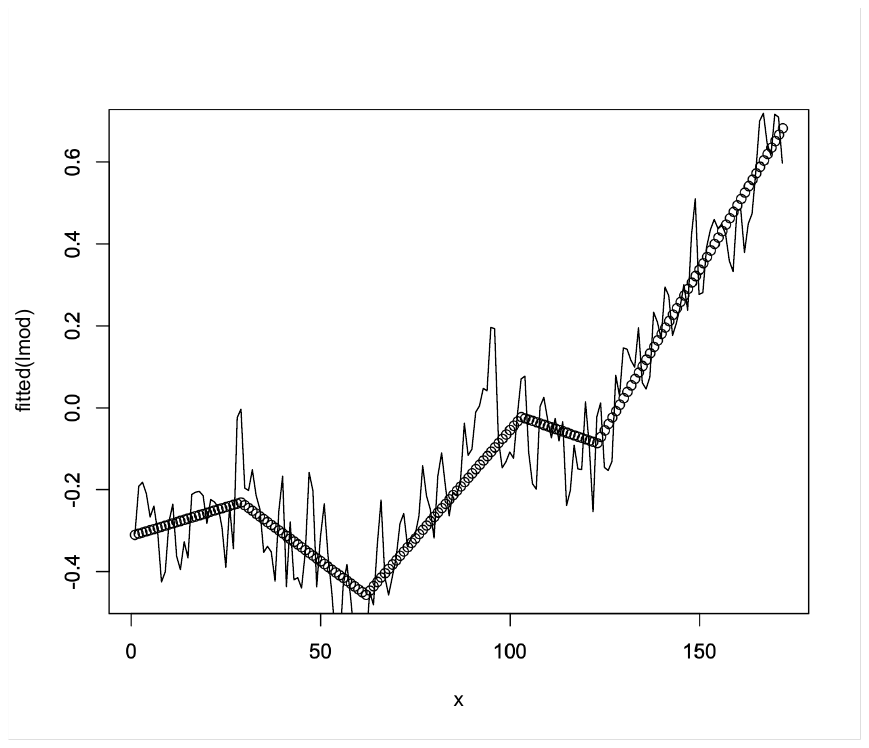} \\
\caption{Linear Analysis of Global Ocean Temperature Anomalies. \label{Figure 5}}
\end{center}
\end{figure}

\medskip\noindent
{\sl Example 7:  AMO and PDO.}
These statistics, Atlantic Multidecadal Oscillation and Pacific Decadal
Oscillation,  monitor
changes in ocean currents in the Atlantic and Pacific Oceans, respectively.  

There are instrumental data for the AMOC 
beginning in 1856 and the PDO beginning in 1900.
In both cases the data have been
recorded monthly, but we average the monthly data to get annual time series
with 163 and 118 values,
respectively.  Since the changes here occur over short time intervals,
with relatively flat stretches between changes, we 
use a model for jump changes.
The methods are those described in our previous paper, although now we allow 
non-zero slopes between change-points and first order
autocorrelation.  

The entire AMO time series produces an estimated value $\hat{\rho} = 0.68$.
All methods agree that there
is a decrease in the mean in 1901, an increase at
1925, a decrease in about  1963, and the latest increase 1995.
A linear analysis incorporating these change-points returns an $R^2$ of about 0.72 and
an auto-correlation of 0.03, which does not test as different from 0.
It seems interesting to observe that the changes in slope detected in the
Northern Hemisphere ocean temperature anomalies discussed in the
preceding example
fall very close to the  mid-points of these intervals.

\begin{figure}
\begin{center}
\includegraphics[width = 5.0 in]{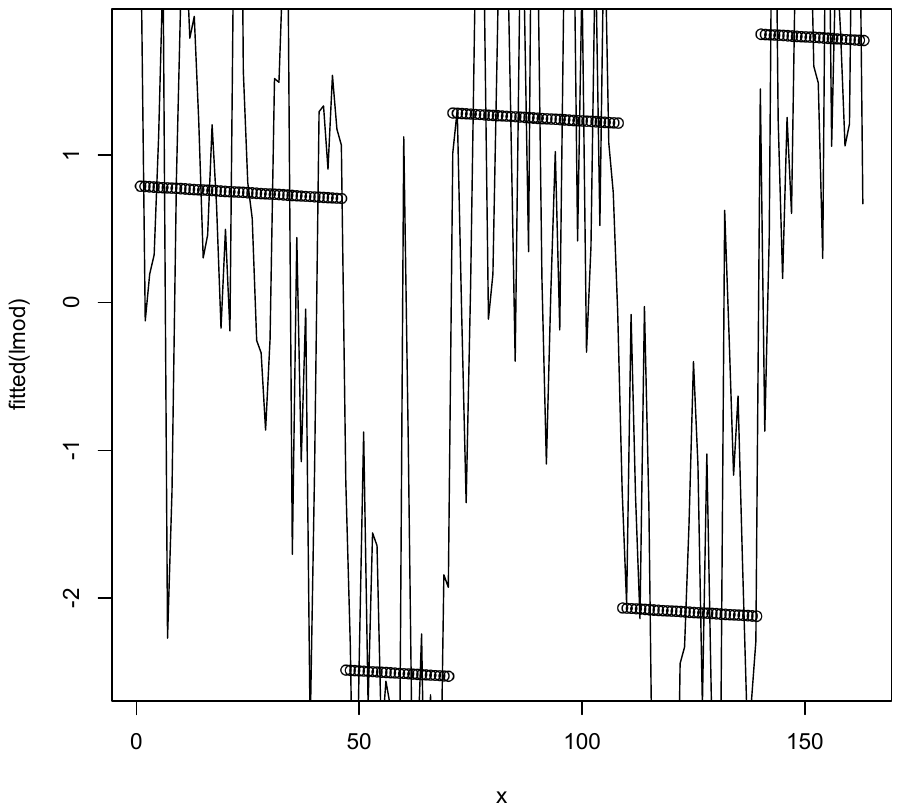} \\
\caption{ Linear Analysis of AMO Data, with Jump Changes of Level in 1901,1925,1963,1995.
\label{Fig 6}}
\end{center}
\end{figure}

The PDO data are similar, but changes appear to be more frequent, as the
name suggests,
and a cumulative sum plot appears very noisy.
Here, if we estimate $\rho$ by
all the observations, we get 0.55 and detect no changes.  If we use
the observations from 1 to 45, to estimate $\rho$, we get  0.30, and we detect changes
in 1947, 1975, 1998, and 2013 by two methods, while
the stricter method that uses all possible background intervals
detects changes in 1947, 1975, and 2013.
A linear analysis with all four detected
change-points produces an $R^2$ of 0.48 and an estimate for $\rho$ of 0.27.

\medskip\noindent
{\sl Example 8: Atmospheric Greenhouse Gases.}
In view of their connection with global warming, it is 
interesting to consider
atmospheric greenhouse gases.  Systematic atmospheric measurements began relatively
recently, although reconstructions from surrogate measurements go back
centuries.  
To illustrate our methods,  we first 
consider global average
methane emissions
for 35 years beginning in 1984,
reported at
ftp://aftp.cmdl.noaa.gov/products/trends/ch4/.

A linear model without changes shows a
positive slope  with an
$R^2$ of 0.95 and an autocorrelation of zero.
MS detects a slowing of the rate of increase in 1993 followed
by a large increase in 2010.  Incorporating these changes into a 
linear model produces an
$R^2$ greater than 0.99 and an autocorrelation of zero.
The statistic (\ref{2locusstatist}) detects the same two changes, 
although it
places the initial decrease about 1996.

Atmospheric measurements of CO$_2$ begin about 1960 and provide a 
similar picture.  Next we consider 514 years, beginning in 1501, of 
surrogate measurements of CO$_2$ provided by 
the Institute for Atmospheric and Climate Science (IAC) at the  
ETH-Z\"urich on the web site
www.co2.earth/historical-co2-datasets.
Unlike other examples, where changes may be missed  
because of large variability,
the variability in these data (and in the recent atmospheric 
measurements) is very
small. There appear to be many small changes that do not help 
us to understand large-scale patterns.
The statistic MS detects 10 changes,  
the first  in 1575 and the last in
1990.  Two changes decrease the slope, while eight increase it.  
See Figure 8.  BIC has a slight preference for a nine change model. 
The top down statistic (\ref{2locusstatist}), which 
detects essentially the same slope changes,  may have some advantages
in problems like this one, since it focuses first on large changes;
and we can stop looking when we think that adding more changes does
not alter our basic insights.  Seq does not work well, since it picks
out many very small changes that it seems reasonable to ignore. 

\begin{figure}
\begin{center}
\includegraphics[width = 5.0 in]{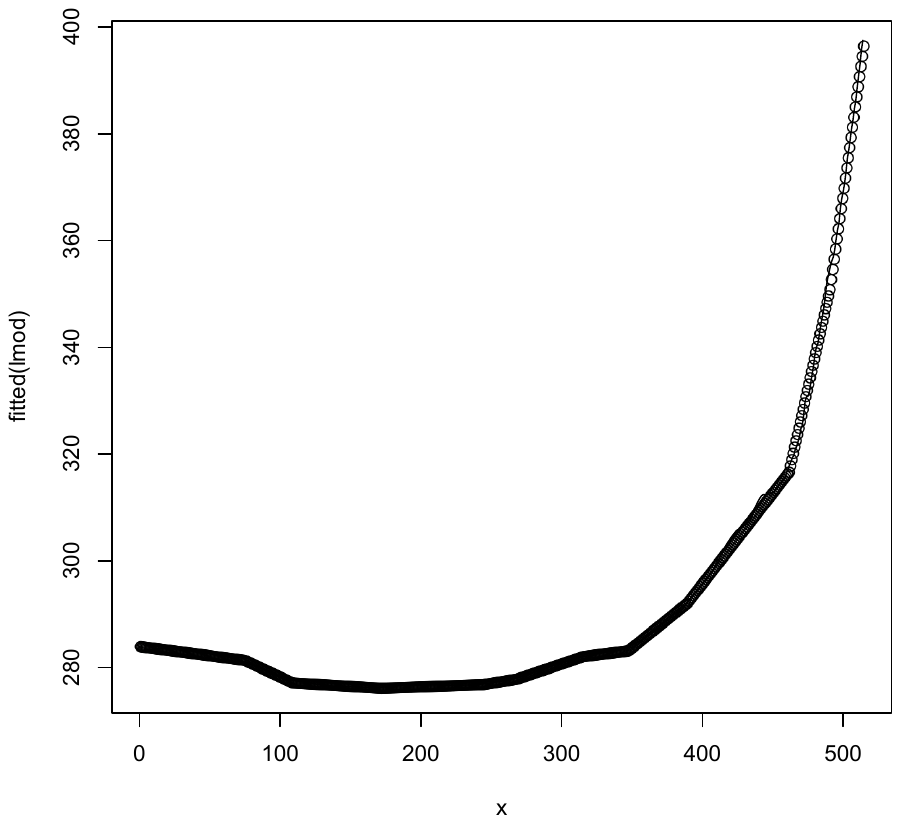} \\
\caption{ Linear Analysis of reconstructed atmospheric CO2, 1501-2014.
\label{Fig 7}}
\end{center}
\end{figure}


\medskip\noindent
{\sl Example 9.} Incidents of hate crimes are reported by the United States FBI web site 
www.fbi.gov from 1996 to 2020.  The record is divided into hate crimes against specific,
overlapping groups, ethnic, relgious, etc.  The largest target of hate crimes is
African-Americans.  A plot showing slope changes in 2008, 2014, and 2019, detected
by (\ref{2locusstatist}) with $\rho = 0$, is given in
Figure 9, where $R^2 = 0.94$, and $\rho$ is estimated to be -0.18.  Seq detects only the
second slope change.  BIC prefers the three change model to a one change model and all 
two change models.  The time series for crimes against Asian-Americans is similar, but the
numbers are smaller, so detection more difficult.  The time series for all hate crimes is similarly
V-shaped, but without the disproportionately large increase between 2019 and 2020.  

\begin{figure}
\begin{center}
\includegraphics[width = 5.0 in]{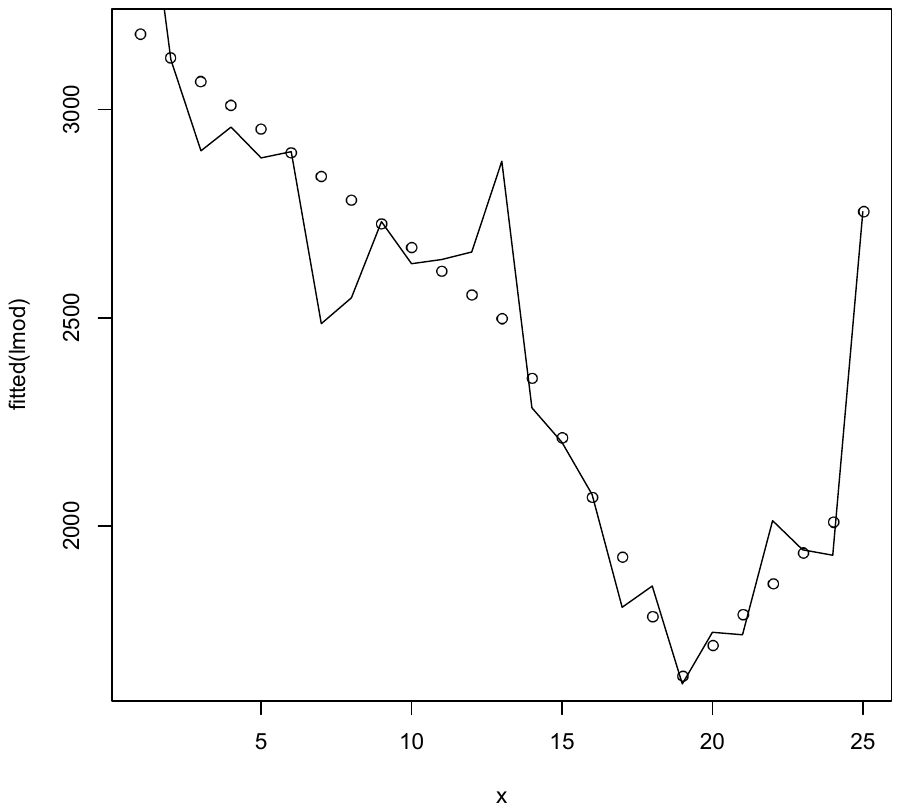} \\
\caption{ Linear Analysis of Hate Crimes against African-Americans.
\label{Fig 8}}
\end{center}
\end{figure}

\bigskip
\centerline{\bf Appendix C: Bump Hunting: fixed shape, variable amplitude and scale.}
In a variety of problems,
the local signal to be detected
occurs in a departure from, followed by a return to, a baseline value.
In many cases the shape of the local signal arises naturally from the
scientific context.
An important example is inherited copy number variation, where there is an
abrupt increase or decrease from the baseline value of two, which is followed
quickly by a return to the
baseline (e.g., \cite{ZhLJS2010}).
Other examples are ChIP-Seq, where a shape to
expect is roughly triangular or double exponential (e.g., \cite{Liu13}) ,
or differential methylation, where a normal probability
density function might be reasonable (cf. \cite{Irizarry12}).

We modify (\ref{loglik}) to include a scale paramter $\tau$,
to obtain the log likelihood function
\begin{equation} \label{loglik2}
-.5\sum_1^T [Y_u - \rho Y_{u-1} -\mu_u -  \sum_k\xi_k f\{(u-t_k)/\tau_k\}]^2.
\end{equation}
Here $f$ is a positive, symmetric integrable function, e.g.,
the square root of (i) a standard normal
density function, (ii) a triangular probability density on $[-1,1]$, (iii)  a
double exponential probability density, or (iv) a uniform probability density
function on $[-1/2, 1/2]$.  Consider the case of at most one
bump and assume that the search interval is long enough relative to the scale
parameter that it is reasonable to ignore end effects.  We
define our standardized
statistic as above, but we now write it as $Z_{t,\tau}$ to reflect the unknown
location $t$ and width $\tau$ of the bump.  An approximation to its false positive
error probability is given in display (9). 

If bumps are sparse, estimation of $\rho$ and $\sigma^2$ is relatively easy,
since occasional signals do not seriously bias maximum likelihood estimators. 

\medskip\noindent
{\sl Remarks.} (i)  The approximation  (9) does not apply to case (iv) mentioned above.
Because of
discontinuities in $f$,  a different approximation
is required (cf. \cite{FaLSi20}).  
(ii) A serious study of these methods is warranted by the variety of
potential applications, some of which involve spatio-temporal 
random fields.  (iii)  For some applications it may be reasonable to assume that
the underlying process is a Poisson or other point process, which among other effects
makes the tail probability of the maximum value of $Z_t$ substantially larger, as 
mentioned in Section 4.  

\centerline{\bf Appendix D:  Additional COVID-19 Examples.}

\noindent{\sl Example 10. COVID-19 in Italy.}
An interesting example of the corona virus is Italy, which had substantial difficulty in
controlling the spread of the virus in the early days of the pandemic.
Using data for 621 days up to (including) 13 October  2021, we find that
(\ref{2locusstatist})  with $\rho = 0.3$ detects
ten slope changes.
A linear least squares analysis produces an $R^2 \approx 0.98$ and an
estimate of $\rho$ of 0.67.   
See Figure 10.  We might want to repeat the analysis with this larger value of $\rho$,
but the effect seems to be to make the change at observation 356 into a borderline case,
which might be omitted,
and leave other slope changes essentially as originally detected.

\begin{figure}
\begin{center}
\includegraphics[width = 5.0 in]{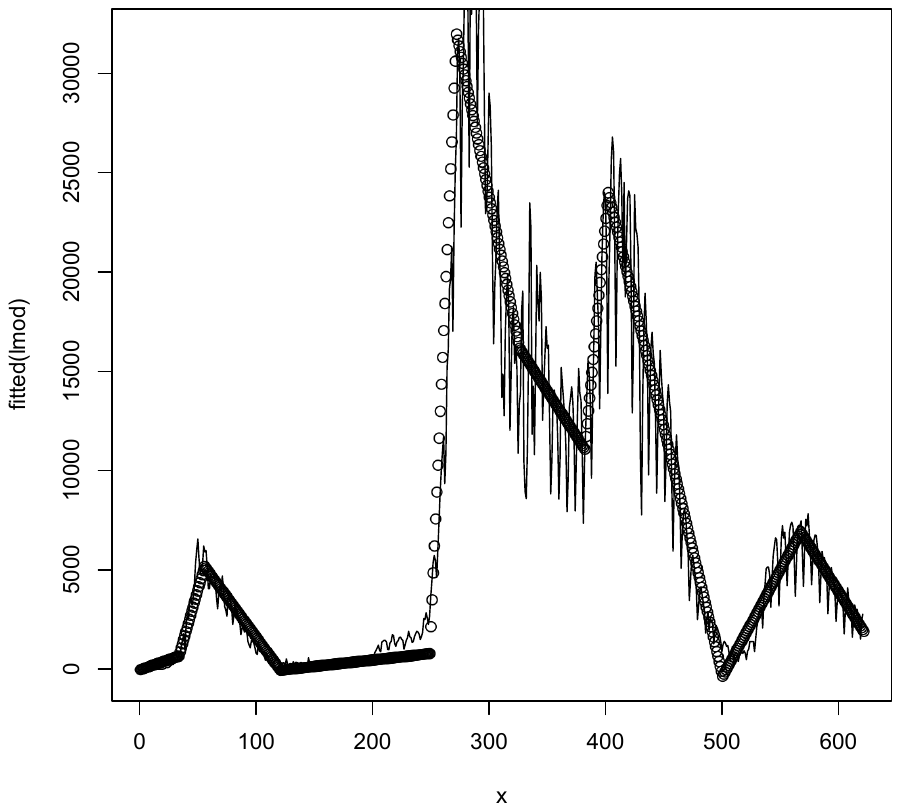} \\
\caption{ Linear Analysis of COVID-19 Incidence in Italy, 31/01/20-25/08/21.
\label{Fig 9}}
\end{center}
\end{figure}

\noindent{\sl Example 11.  COVID-19 in Hong Kong}
An extremely well organized and informative web site is 

\centerline{chp-dashboard.geodata.gov.hk/covid-19/en.html};

Ourworldindata also contains well organized daily incidence data for Hong Kong, although it did 
not in the early days of the pandemic.

\smallskip\noindent
For 579 days from 25/01/2020 to 25/08/2021 all three methods give about the same
results.  MS produces the best looking broken line, largest $R^2$ 
(about 0.94) and smallest $\rho$ (about 0.35) in the linear analysis.  
See Figure  11.

\begin{figure}
\begin{center}
\includegraphics[width = 5.0 in]{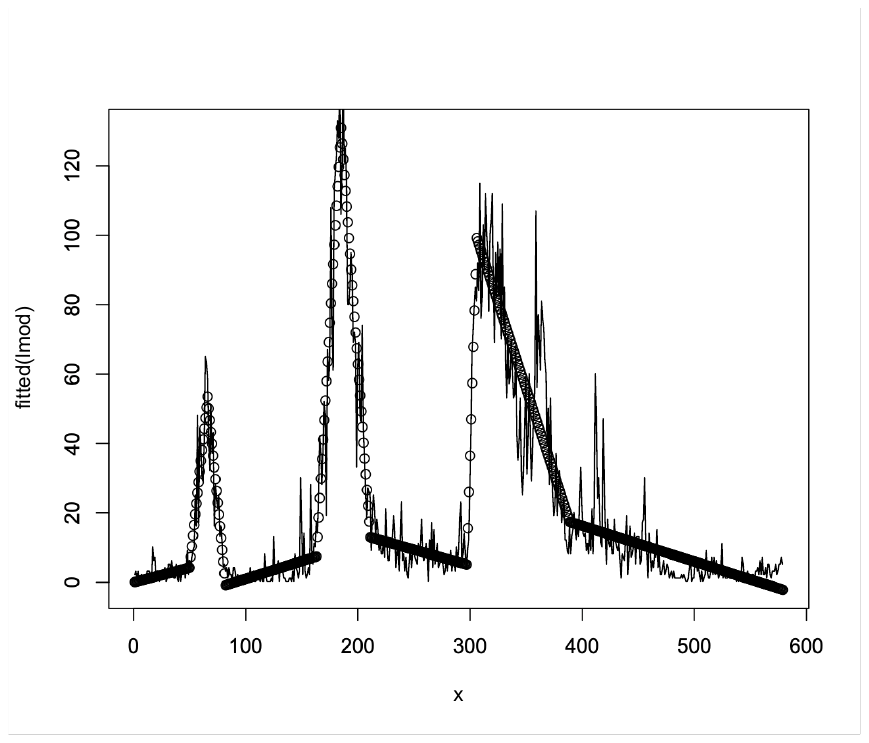} \\
\caption{ Linear Analysis of COVID-19 Incidence in Hong Kong, 25/01/20-25/08/21.
\label{Fig 10}}
\end{center}
\end{figure}

For daily incidence of COVID-19, it is apparent from Figures 10 and 11 that the variability is much larger
when the incidence is large.  It does not seem important to deal with this heteroscedasticity because
slope changes are also usually much larger when the incidence is large.  We have, however,  considered 
the possibilities of transforming the data, say by taking a square root or a logarithm.  This rarely
seems to have a consequential effect.  

Since simple epidemic models suggest that the incidence of an infections disease  
increases (and decreases) exponentially, one might prefer using a (usually overdispersed) log Poisson 
linear model.  This sometimes produces a visually more appealing plot and in some cases allows one
to discard some change-points because the natural curvature of the model substitutes for a slope
change in a broken line model.

\bigskip
\centerline{\bf Appendix E:  Threshold Autoregression}
As a final illustration of our methods we consider a simple case of Threshold Autoregression,
which has been studied in a number
of papers by Tong and colleagues.

For a model assume
\begin{equation} Y_u = \mu + \rho Y_{u-1} + \xi Y_{u-1} I\{ Y_{u - 1} \leq t \} + \epsilon_u
\end{equation}

To test $\xi = 0$, in the notation of Section 2
\[
\ell_\xi(t) = \sum_u (Y_u - \mu -\rho Y_{u-1}) Y_{u-1} I\{ Y_{u-1} \leq t\},
\]
and asymptotically under the null hypothesis we have
\[\Psi(t)' \sim
T (\Expec(Y;Y \leq t),\Expec(Y^2;Y \leq t)),
\]
\[
\dot{\ell}_\theta = (\sum (Y_u - \mu - \rho Y_{u-1}), \sum [(Y_u - \mu - \rho Y_{u-1}) Y_{u-1}],
\] and $A^{-1}$ is
the $2\times 2$ matrix with entries $a_{11} = T, \; a_{12} = T \Expec(Y), \; a_{22} = T \Expec(Y^2),$
where $Y$ denotes the stationary
distribution of $Y_u$ under the null hypothesis (and the stationarity
assumption that $|\rho| < 1$).
Hence
\[\sigma^2(t) \sim T \Expec(Y^2;Y \leq t)
\]
\begin{equation} \label{var}
 - T[ \Expec^2(Y^2; Y\leq t) +
\Expec(Y^2) \Expec^2(Y;Y \leq t) - 2 \Expec(Y) \Expec(Y;Y \leq t) \Expec(Y^2;Y \leq t)]/\Var(Y).
\end{equation}

Let $G(t) =
\Expec(Y^2;Y \leq t)$, so
$\sigma^2(t) \sim T [G(t) - \Psi(t)'A\Psi(t)].$  Straightforward
calculations show that
the covariance of the standardized score statistics $Z_s$ and $Z_t$
are given by $[G(\min(s,t)) - \Psi'(s)A\Psi(t)]/\sigma(s) \sigma(t).$
Hence locally for small $\delta$

\[ \Cov(Z_t, Z_{t+\delta}) \approx I-(\delta/2)[G - \Psi' A \Psi]^{-1}
\dot{G},
\]
where the functions $G$, $\Psi$  and $\dot{G}$ are evaluated at $t$,
and slightly more generally
\[\Cov[(Z_{t+\delta_1}, Z_{t+\delta_2})|Z_t]  \approx \min(\delta_1, \delta_2)
[G - \Psi' A \Psi]^{-1}\dot{G}.
\]
This allows us to calculate an approximation to
$\Prob\{\max_t |Z_t| \geq b\}$.
The methods of, e.g., Woodroofe (1976) or Yakir (2013) lead after some calculation to
\begin{equation} \label{tarapprox}  \Prob \{\max_{t_0 < t < t_1} |Z_{t}| > b \}
\approx 2 b \varphi(b)\int_{t_0}^{t_1} [ \Psi'(t)A\dot{\Psi}(t)/\sigma^2(t)] dt.
\end{equation}

We have evaluated (\ref{tarapprox}) by numerical integration, and the
p-values given below 
come from that evaluation with the observed mean value
and standard deviation.  We have also implemented an empirical version of this
approximation, where the entire computation is based on the appropriately
estimated quantities.  These two approximations are in rough agreement for
our examples, where the p-value is very small.
Simulations suggest that the Type I error control of (\ref{tarapprox})
is adequate under the model
The power also seems reasonable.  In principle, one should learn something
from  the maximizing point of the statistic $Z_s$  relative to the estimated values
of $\mu$ and of $\sigma$, but in simulations the maximizing value of shows more
variability than we can easily interpret.
A calculation similar to that giving (\ref{var}) indicates that in the case that
$\xi \neq 0$
the expected value of $Z_t$ is approximately $\xi T^{1/2} \sigma(t)$, as it would be for
likelihood theory with standard regularity conditions.  This approximation seems relatively
stable and suggests a rough approximation for $\xi$ by equating the approximation to
the observed value of $\max_s Z_s$ (or $\min_s Z_s$).

With a slight modification in principle, but substantially more
detailed calculation in application, one can consider higher
order autoregressions.  Suppose, for example, the model is
\[
Y_u = \rho_1 Y_{u-1} + \rho_2 Y_{u-2} + (\xi_1 Y_{u-1}
+\xi_2 Y_{u-2})I\{Y_{u-2} \leq t\}  + \epsilon_u.
\]
Now $X_t$ is a two-dimensional vector, while $G(t)$, $\Psi(t)$,
and $\sigma^2(t) = \Sigma_t,$ say, are
$2\times2$ matrices.
An appropriate test statistic has the form
$\max_t V_t' \Sigma^{-1}_t V_t.$  Putting
$Z_t = \Sigma^{-1/2} V_t$, we can control the false positive
probability based on the approximation
\[\Prob \{\max_{t_0 < t < t_1} ||Z_{t}|| > b \}
\]
\begin{equation} \label{approxthreshold2}
\approx (b^2/2) \exp(-b^2/2) (2 \pi)^{-1} \int_{t_0}^{t_1} \int_0^{2 \pi}  {\bf e}'
\Sigma^{-1}_t \dot{G}(t)  {\bf e}
d \omega dt + \exp(-b^2/2),
\end{equation}
where ${\bf e}  = (\cos(\omega), \sin(\omega))'$, and where
$\pi^{-1}$ times the integral over $\omega$ equals
the trace of $\Sigma^{-1}_t \dot{G}(t)$.

Simulations suggest that the Type I error control of this
procedure is adequate
under the model.
The power also seems reasonable.  In principle, one should learn something
from  the maximizing point of the statistic $Z_t$  relative to the estimated values
of $\mu$ and of $\sigma$, but in simulations the maximizing value of $t$ shows more
variability than we can easily interpret.

\medskip\noindent
{\sl Examples.}
As illustrative applications, we follow Chan and Tong (1990) in considering the
Canadian lynx data and
Nicholson's blowfly data.  We do not, however, try to give a
complete discussion of these well studied data.

The Canadian lynx data consist of annual counts over a
114 year period, which are considered a surrogate for the size of the lynx population.
The scientific reasoning behind a two-phase model is the hypothesis
that when the lynx
population is small, it finds a more than adequate food supply and increases
until it reaches a
size where the food supply is inadequate, when it then decreases
until the cycle
begins again.  Empirical observations suggest that the decrease in times of food
shortage is more rapid than the increase in times of abundance.
If we use a first order auto regressive model, these fluctuations suggest
a positive autoregressive parameter in times of food abundance, which
decreases when there is a food shortage.
Although there are clear outliers in the data, which has lead others to
consider transformations to shorten the tails, we analyze
the original data.
Employing the first order autoregressive model with a
possible change in the autoregressive parameter, as suggested above
yields null estimators of $\Expec(Y)$ and $\Var(Y)$ of
1538 and 1560, respectively.
The maximum value of $Z_t$ is 3.89, which
occurs at about $t = 3500$, and the resulting (two-sided) p-value
based on (\ref{tarapprox}) is
approximately $0.004$.

Nicholson's blowfly data (\cite{Brill80} from the web site

\smallskip\noindent
\centerline{www.stat.berkeley.edu/~brill/blowfly97I.html},

\smallskip\noindent
contains 361 observations in four columns, labeled respectively ``births,''
``nonemerging,'' ``emerging,'' and ``deaths.''  The first, third,
and  fourth columns vary
from numbers close to 0 to numbers in the several thousands.
The second column typically involves substantially smaller numbers.
For the first column.
the estimated mean value, autocorrelation,  and standard deviation
under the null model are 1438, 0.6 and 1670, respectively.
The maximum $Z$ value  is about 4.73, which occurs for $t \approx 5300.$
The p-value for the hypothesis of no change
is approximately $2 \times 10^{-4}.$  For ``emerging''
the maximum $Z$-value is about 5.22, which gives a p-value of
about $6 \times 10^{-5}.$  For these data the change in the autocorrelation is
positive.

\bibliography{Bibliography-MM-MC}

\begin{thebibliography}{}

\bibitem[Adler and Taylor(2007)]{AT07}
R.~J. Adler and J. ~E. Taylor (2007).
\newblock \textit{Random Fields and Geometry}
\newblock Springer-Verlag, New York-Heidelberg-Berlin.

\bibitem[Arrhenius(1896)]{Ar96}
S. Arrhenius (1896).  
\newblock On the influence of carbonic acid in the air
on the temperature on the ground, 
\newblock \textit{Philos. Mag.} \textbf{41}, 237-276.

\bibitem[Baranowski, Chen, and Fryzlewicz(2019)]{BCFr16}
R. ~Baranowsk, Yining ~Chen, and P. ~Fryzlewicz (2019).
\newblock Narrowest-over-threshold detection of multiple change-points
and change-point-like features
\newblock \textit{J. Roy. Statist. Soc. B} \textbf{89},
649-672.



\bibitem[Chan and Tong (1990)]{CT90}
K. ~S. Chan and H. Tong (1990)
\newblock On likelihood ratio tests for threshold autoregression
\newblock \textit{JRSSB} \textbf{52}, 469-476.


\bibitem[Davies(1987)]{Dav87}
R. ~B Davies (1987).  
\newblock Hypothesis testing when a nuisance parameter is
present only under the alternative 
\newblock \textit{Biometrika} \textbf{ 74}, 33-43.

\bibitem[Fang, Li and Siegmund(2020)]{FaLSi20}
X.~Fang, J. ~Li  and D.~O. Siegmund (2020).
\newblock Segmentation and estimation of change-point models:
false positive control and confidence regions
\newblock \textit{Ann. Statist.} \textbf{43}, 1615--1647.

\bibitem[Fryzlewicz(2014)]{Fr14}
P.~Fryzlewicz (2014).
\newblock Wild binary segmentation for multiple change-pont detection.
\newblock \emph{Ann. Statist.} \textbf{42},  2243--2281.

\bibitem[Heyde, 1997]{Heyde97}
C.~C. Heyde (1997)
\newblock \textit{Quasi-Likelihood and Its Application}
\newblock Springer-Verlag, New York-Heidelberg-Berlin.

\bibitem[Jaffe, {\sl et al.} (2012a) ] {Irizarry12}
A. ~E. Jaffe, A. ~P. Feinberg, R. ~A. Irizarry, and J ~T. Leek (2012a).
\newblock Significance analysis and statistical dissection of
variably methylated regions.
\newblock \textit{Biostatistics} \textbf{13}, 166--178.

\bibitem[Jaffe, {\sl et al.} (2012b) ] {JaffeIrizarry12}
A. ~E. Jaffe,  P. Murakami,  Hwajin Lee, J ~T. Leek, 
M. Daniele Fallin, A ~P. Feinberg, and R. ~A. Irizarry (2012b).
\newblock Bump hunting to identify differentially methylated regions
in epigenetic epidemiology studies,
\newblock \textit{Int. J. of Epidemiology} \textbf{41}, 200--209. 

\bibitem[Jones, {\sl et al.} (2012)] {Jonesetal12}
P. ~D. Jones, D. ~H. Lister, T. ~J. Osborn, C. Harpham, M. Salmon, and C. ~P. Morice
\newblock Hemispheric and large-scale land surface air temperature variations:
An extesive revision and an update to 2010.
\newblock {\sl J. Geophys. Res.} \textbf{117} D05127,doi:10.1029/2011JD017139.

\bibitem[Knowles and Siegmund (1989)]{KnSi89}
M. Knowles and D. Siegmund (1989)
\newblock On Hotelling's approach to testing for a nonlinear parameter in regression
\newblock {\sl Int. Statist. Rev.} \textbf{57} 205-220.

\bibitem[Knowles, Siegmund and Zhang (1991)]{KnSZh91}
M. Knowles, D. Siegmund and H. ~P. Zhang (1991)
\newblock  Confidence regions in semilinear regression 
\newblock {\sl  Biometrika}  \textbf{79} 15-31.

\bibitem[Krieger, Pollak and Yakir (2003)]{KPY03}
A. ~B. Krieger, M. Pollak, and B. Yakir (2003)
\newblock Surveillance of a simple linear regression
\newblock J. Amner. Statist. Assoc. \textbf{98} 456-469.


\bibitem[Manley (1974)]{Man74}
Gordon Manley (1974)
\newblock Central England temperatures:monthly means 1659 to 1973
\newblock \textit{Quarterly J. Roy. Meteorological Soc.} \textbf{100} 389-405.

\bibitem[Muggeo (2016)] {Mug16}
V. ~M. ~R Muggeo (2016)
\newblock Testing with a nuisance parameter present only under the alternative:
a score based approach with application to segmented modeling
\newblock \textit{J. Statistical Computation and Simulation} doi: 10.1080/00949655.2016.1149855

\bibitem[Olshen {\sl et al.}(2004)]{OlVe04}
A.~B. Olshen, E.~S. Venkatraman, R.~Lucito and M.~Wigler (2004).
\newblock Circular binary segmentation for the analysis of
array-based DNA copy number data.
\newblock \textit{Biostatistics} \textbf{5}, 557--572.

\bibitem[Raftery and Akman (1986)]{RafAk86}
A. ~E. Raftery and V. ~E. Akman (1986).
\newblock Bayesian analysis of a Poisson process with a change-point.
\newblock \emph{Biometrika} \textbf{73}, 85--80.

\bibitem[Robbins, Gallagher, and Lund  (2016)]{Lund16}
M. ~W. Robbins, C. ~M. Gallagher, and R. Lund  (2016).
\newblock A general regression change-point test for time series data,
\newblock \emph{Jour. Amer. Statist. Assoc.} {\bf 111} 670--683.

\bibitem[Schwartzman \; {\sl et al.} (2013)]{SJGM2013}
A. Schwartzman, A. Jaffe, Y.  Gavrilov, Y. and  C. ~ E. Meyer (2013).
\newblock Multiple testing of local maxima for 
detection of peaks in ChIP-seq data,
\newblock \emph{Ann. Appl. Statist.}, \textbf{7}, 471-494.

\bibitem[Shin \; {\sl et al.} (2013)] {Liu13}
H. ~Shin, T. ~Liu, X. ~Duan, Y. ~Zhang, X. ~S. Liu (2013).
\newblock Computational methodology for ChIP-seq analysis.
\newblock \emph{Quantitative Biology} \textbf{1}, 54--70.

\bibitem[Siegmund and Yakir(2000)]{SY2000}
D.~O. Siegmund and B. Yakir (2000).
\newblock Tail probabilities for the null distribution
of scanning statistics.
\newblock
\emph{ Bernoulli}  \textbf{6}, 191--213.

\bibitem[Siegmund and Worsley(1995)]{SiWo95}
D.~O. Siegmund and K.~J. Worsley (1995).
\newblock Testing for a signal with unknown location and scale 
in a stationary Gaussian random field. 
\newblock \emph{Ann. Statist.} \textbf{23}, 608--639.

\bibitem[Smith and Cook (1980)]{SC80}
A. ~F. ~M. Smith and D. ~G. Cook (1980).
\newblock Strait lines with a change-point:  an analysis of some renal transplant data
\newblock {\sl J. Roy. Statist. Soc. Series C} \textbf{29}, 180-189.

\bibitem[Toms and Lesperance (2003)]{ToLes03}
J. ~ D. Toms and M. ~L. Lesperance
\newblock Piecewise regression:  a tool for identifying ecological thresholds
\newblock {\sl Ecology} \textbf{84}, 2034-2041.

\bibitem[Zhang {\sl et al.}(2010)] {ZhLJS2010}
N.~R. Zhang, D.~O. Siegmund, H.~Ji and J.~Z. Li (2010).
\newblock Detecting simultaneous changepoints in multiple
sequences. With supplementary data available online.
\newblock \emph{Biometrika} \textbf{97}, 631--645.

\end{thebibliography}

\begin{thebibliography}{}
\bibitem[Brillinger, et al. 1980]{Brill80}
D. ~R. Brillinger, J. Guckenheimer, P. Guttorp and G. Oster (1980)
\newblock Empirical modelling of population time series data: the case of age and density dependent vital rates
\newblock \textit{Lectures on Mathematics in the Life Sciences} \textbf{13}, 65-90
\newblock American Math. Soc.

\bibitem[Chan (1991)]{KSC91}
K. ~S. Chan (1991)
\newblock Percentage points of the likelihood ratio test for threshold autoregression
\newblock \textit{JRSSB} \textbf{53}, 691-696.

\bibitem[Chan and Tong (1990)]{CT90}
K. ~S. Chan and H. Tong (1990)
\newblock On likelihood ratio tests for threshold autoregression
\newblock \textit{JRSSB} \textbf{52}, 469-476.

\bibitem[Church and White(2011)]{CW11}
J. ~A. Church and N. ~J. White (2011)
\newblock Sea level rise from the late 19th to the early 21st century
\newblock \textit{ Surv. Geophysis.} \textbf{32}, 585-602.

\bibitem[Jaffe, {\sl et al.} (2012a) ] {Irizarry12}
A. ~E. Jaffe, A. ~P. Feinberg, R. ~A. Irizarry, and J ~T. Leek (2012a).
\newblock Significance analysis and statistical dissection of
variably methylated regions.
\newblock \textit{Biostatistics} \textbf{13}, 166--178.

\bibitem[Jaffe, {\sl et al.} (2012b) ] {JaffeIrizarry12}
A. ~E. Jaffe,  P. Murakami,  Hwajin Lee, J ~T. Leek,
M. Daniele Fallin, A ~P. Feinberg, and R. ~A. Irizarry (2012b).
\newblock Bump hunting to identify differentially methylated regions
in epigenetic epidemiology studies,
\newblock \textit{Int. J. of Epidemiology} \textbf{41}, 200--209.
 
\end{thebibliography}

\end{document}